\documentclass[a4paper]{amsart}

\usepackage{amsmath, amsfonts, amsthm, framed}
\usepackage{xspace}
\usepackage{newunicodechar}
\usepackage{fullpage}
\usepackage{enumitem}
\usepackage{quiver}
\usepackage{fontawesome}
\usepackage[numbers]{natbib}
\usepackage{soul}
\usepackage{cancel}
\usepackage{xspace}
\usepackage{eucal}
\usepackage{colonequals}
\usepackage[colorlinks=true]{hyperref}
\definecolor{darkergreen}{rgb}{0.0, 0.5, 0.0}
\definecolor{Blue}{RGB}{0,0,148}

\usepackage{listings}

\definecolor{light-gray}{gray}{0.95}
\usepackage{tikz}
\usetikzlibrary{positioning, arrows.meta, shapes.geometric}
\usepackage[a4paper, margin=0.75in]{geometry} % More page space

\definecolor{keywordcolor}{rgb}{0.7, 0.1, 0.1}   % red
\definecolor{commentcolor}{rgb}{0.4, 0.4, 0.4}   % grey
\definecolor{symbolcolor}{rgb}{0, 0, 0.8}    % blue
\definecolor{tacticcolor}{rgb}{0, 0, 0.8}    % blue
\definecolor{sortcolor}{rgb}{0.1, 0.5, 0.1}      % green
\newcommand*{\lean}[1]{\lstinline{#1}\xspace} % Lean in-lineexpressions
\theoremstyle{plain}
\newtheorem{theorem}{Theorem}[section]

\newtheorem{lemma}[theorem]{Lemma}
\newtheorem{proposition}[theorem]{Proposition}
\newtheorem{corollary}[theorem]{Corollary}
\theoremstyle{definition}
\newtheorem{definition}[theorem]{Definition}

\newcommand{\height}{\mathrm{ht}}
\newcommand{\Ann}{\mathrm{Ann}}
\newcommand{\Ass}{\mathrm{Ass}}
\newcommand{\Ker}{\mathrm{ker}}

\newcommand{\Supp}{\mathrm{Supp}}
\DeclareMathOperator{\Spec}{Spec}
\DeclareMathOperator{\projdim}{projdim}
\DeclareMathOperator{\injdim}{injdim}
\DeclareMathOperator{\gldim}{gldim}
\DeclareMathOperator{\depth}{depth}
\DeclareMathOperator{\Ext}{Ext}
\DeclareMathOperator{\Tor}{Tor}
\newcommand*{\mathlib}{\textsc{mathlib}\xspace} %% Temporarily set this way to catch forgotten instances
\newcommand*{\Lean}{\textsc{Lean}\xspace} %% Temporarily set this way to catch forgotten instances

\newcommand*{\Leanf}{\textsc{Lean4}\xspace}
\title{Formalization of Auslander--Buchsbaum--Serre criterion in Lean4}

\author{Nailin Guan}
\email{nailinguan55@stu.pku.edu.cn}
\author{Yongle Hu}
\email{huyongle@stu.pku.edu.cn}

\begin{document}

\begin{abstract}
We present a comprehensive formalization in the Lean4 theorem prover of the Auslander--Buchsbaum--Serre criterion, which characterizes regular local rings as those Noetherian local rings with finite global dimension. Rather than following the well-known proof that computes the projective dimension of the residue field via quotient by regular sequences and uses the Koszul complex to bound the cotangent space dimension by the global dimension, our approach is built systematically on the formalization of depth defined via the vanishing of Ext functors. We establish key homological results including Rees' theorem, the Auslander--Buchsbaum formula, and Ischebeck's theorem, and further develop the theories of Cohen--Macaulay modules and rings, including a complete formalization of the unmixedness theorem for Cohen--Macaulay rings. To prove the Auslander--Buchsbaum--Serre criterion, we show that maximal Cohen--Macaulay modules over regular local rings are free and establish a weakened form of the Ferrand--Vasconcelos theorem specific for the unique maximal ideal. As corollaries, we deduce that regularity can be checked at maximal ideals and formalize Hilbert's Syzygy Theorem. This work demonstrates how homological algebra can be effectively employed in the formalization of commutative algebra, providing extensive infrastructure for future developments in the field.
\end{abstract}

\maketitle

\section{Introduction}

The Auslander--Buchsbaum--Serre criterion is a beautiful mathematical result stating that for a Noetherian local ring $R$, the ring $R$ is regular if and only if it has finite global dimension. (In this paper, all rings are assumed to be commutative.) This criterion was not proved all at once. At first, only the “only if” direction was known, as shown in \cite{Auslander1953}. At that time, it was still an open problem whether the localization of a regular local ring at a prime ideal remains regular. The homological characterization of regular local rings was completed only after Jean--Pierre Serre proved the converse direction later in \cite{Serre1956}, as formally recorded in \cite{AuslanderBuchsbaum1956}.
Once this characterization became available, it became straightforward to prove that the localization of a regular local ring is again regular. This was one of the earliest examples where homological methods were used to prove a result that is not itself homological. It shows how homological invariants—such as projective dimension, depth, and homological dimension—can detect and control purely “geometric” or “structural” properties of rings. This opened the door for homological algebra to become a central tool in commutative algebra and algebraic geometry. A similar homological perspective later led to the Auslander--Buchsbaum formula, Serre's intersection multiplicity theory \cite{intersection65}, the development of Cohen--Macaulay and Gorenstein ring theory, and deformation theory, for example \cite{deformation87}.

This project began with the goal of developing the theories of depth and Cohen--Macaulay rings. As the concepts of regular elements and regular sequences were already formalized, this goal was within reach using the current library \mathlib{} \cite{mathlib}. The effort was not only of mathematical interest, but also a test of the newly developed homological algebra in \mathlib{} via its application to commutative algebra. As we shifted our focus toward combining it with homological algebra, we set Hilbert's Syzygy theorem and the Auslander--Buchsbaum--Serre criterion as our main targets. However, at the time of this work, the $\Tor$ functor was still largely missing from \mathlib{}, and the development of the Koszul complex was also incomplete, so the most well-known proof could not be used in the formalization. The turning point was finding an approach to the global dimension of a regular local ring by proving that a maximal Cohen--Macaulay module over a regular local ring is free; then the target became attainable. Later, inspired by the proof of the Ferrand--Vasconcelos theorem in \cite[Theorem 2.2.8]{CM_ring}, we realized that induction on the span rank (minimal number of generators) could still be carried out by focusing only on the unique maximal ideal. This allowed us to complete the full Auslander--Buchsbaum--Serre criterion.
The value of this work lies not only in the final theorems: throughout the process, we formalized more basic constructions in commutative and homological algebra (for example, the definitions of projective dimension and $\depth$), tested these definitions in later developments, and improved their interfaces to make them easier to use. However, some of these additional interfaces are not included in this paper due to space limitations.
Currently, all of section \ref{sec:pre} except for subsection \ref{subsec:gldimpre} (about global dimension) is already in \mathlib{}. We plan to upstream most of our work into \mathlib{} in the future.

\subsection*{Outline of the Paper}

The main goal of the formalization project discussed in this paper is the Auslander--Buchsbaum--Serre criterion.

\begin{definition}
A Noetherian local ring $(R, \mathfrak{m})$ is a \emph{regular local ring} if the span rank (minimal number of generators) of $\mathfrak{m}$ is equal to the Krull dimension of $R$.
\end{definition}

Recall that the projective dimension of a module is the minimal length of its projective resolution, and the global dimension of a ring $R$ is the supremum of the projective dimensions of all $R$-modules, denoted by $\gldim(R)$. The criterion can be stated as follows (see also Subsection ~\ref{subsec:concl}).

\begin{theorem}[Auslander--Buchsbaum--Serre]
    A Noetherian local ring $R$ is regular if and only if $\gldim(R)$ is finite.
\end{theorem}

The proof of the Auslander--Buchsbaum--Serre criterion splits into two largely independent parts, corresponding to the two directions of the equivalence.
\begin{theorem}\label{thm:reggldim'}
    For a regular local ring $R$, $\gldim(R) = \dim(R)$.
\end{theorem}
\begin{theorem}\label{thm:regofproj'}
    For a Noetherian local ring $(R, \mathfrak{m}, k)$, if $k$ has finite projective dimension, then $R$ is regular.
\end{theorem}

An important corollary of Theorem \ref{thm:reggldim'} is the following theorem. 
\begin{theorem}[Hilbert's Syzygy]
    For a field $k$, $\gldim(k[x_1, \cdots, x_n]) = n$.
\end{theorem}

%In Section ~\ref{sec:pre}, we introduce the preliminary results needed for the formalization, including existing theories in \mathlib{} and supplementary results we developed during the study of Cohen--Macaulay rings and regular local rings. We then present the development of depth and Cohen--Macaulay modules and rings in Sections ~\ref{sec:depth} and ~\ref{sec:CMdef}, following \cite[Section 6]{mats_commalg_2ed}. This was the earliest part to be formalized. 

As for the proof of the criterion, i.e., Theorem \ref{thm:reggldim'} and Theorem \ref{thm:regofproj'}, the key preliminary ingredient is the Auslander--Buchsbaum formula (Theorem ~\ref{thm:ABthm}) in Subsection ~\ref{subsec:ABthm}. Recall for a (finitely generated) module $M$ over a local ring $R$, $\depth(M)$ is the maximal length of $M$-regular sequences.

\begin{theorem}[Auslander--Buchsbaum]\label{thm:ABthm'}
Let $(R, \mathfrak{m})$ be a Noetherian local ring and let $M$ be a non-zero finitely generated $R$-module. If $\projdim(M)$ is finite, then
$$
\projdim(M) + \depth(M) = \depth(R).
$$
\end{theorem}

\begin{itemize}[leftmargin=*]
\item For the proof of Theorem \ref{thm:reggldim'}, consider a regular local ring $(R, \mathfrak{m}, k)$ of dimension $d$.
\end{itemize}

By Theorem ~\ref{thm:ABthm'}, the projective dimension of any finitely generated module is either infinite or at most $\depth(R)=d$, and since $\depth(k)=0$, this upper bound is attained for $k$. This reduces the problem to showing that all finitely generated modules over a regular local ring have finite projective dimension.
\begin{definition}
    A module $M$ over a local ring $R$ is a \emph{maximal Cohen--Macaulay module} if $\depth(M) = \dim(R)$.
\end{definition}
By basic inequalities of $\depth$, we obtain the following:
\begin{proposition}
    The $d$-th syzygy of a finitely generated module $M$ in a free resolution is either a maximal Cohen--Macaulay module or $0$.
\end{proposition}

Here and throughout this article, the $n$-th syzygy of a free resolution of a module $M$ means kernel of the $(n-1)$-th differential, with the convention that the $0$-th differential is the surjection to $M$.

The following result about maximal Cohen--Macaulay modules then ensures that this yields a finite projective (free) resolution of $M$.
\begin{theorem}[{\cite[00NT]{stacks-project}}]\label{thm:maxfree'}
    For a regular local ring $R$, any finitely generated maximal Cohen--Macaulay module over $R$ is free.
\end{theorem}

\begin{itemize}[leftmargin=*]
\item For the converse direction (Theorem \ref{thm:regofproj'}), i.e. showing regularity from finite global dimension, we initially attempted to use the following theorem.
\end{itemize}

\begin{theorem}[Ferrand--Vasconcelos, {\cite[Theorem 2.2.8]{CM_ring}}]\label{thm:FVthm'}
    For an ideal $I$ of a Noetherian local ring $R$, if $I$ has finite projective dimension and $I/I^2$ is free over $R/I$, then $I$ is generated by a regular sequence.
\end{theorem}

However, since some techniques about calculating rank in a free resolution used in its proof are not yet fully formalized, we focus only on the case where $I$ is the unique maximal ideal, which still suffices for our goal. Note that a regular sequence in $\mathfrak{m}$ has length at most $\depth(R) \leq \dim(R)$, this would force $R$ to be regular.

Returning to the case $I = \mathfrak{m}$, the proof roughly follows a modified version combining ideas from \cite{stacks-project} and \cite[Section 2.2]{CM_ring}. By induction on the minimal number of generators of $\mathfrak{m}$, if $x$ is a regular element in $\mathfrak{m} \setminus \mathfrak{m}^2$, then $\mathfrak{m}/(x)$ is a direct summand of $\mathfrak{m}/x\mathfrak{m}$. The following proposition then implies that $\mathfrak{m}/(x)$ has finite projective dimension, so the induction hypothesis applied to $R/(x)$ completes the proof.
\begin{proposition}
    For a Noetherian local ring $(R, \mathfrak{m})$ and a finitely generated $R$-module $M$ and an element $x \in \mathfrak{m}$ that is both $R$-regular and $M$-regular, we have
    $$\projdim_{R/(x)}(M/xM) = \projdim_R(M).$$
\end{proposition}

\noindent It is worth noting that properties of Cohen--Macaulay rings are established systematically in Section ~\ref{sec:CMdef}, and not all results presented here are directly related to the Auslander--Buchsbaum--Serre criterion. However, most of these unrelated results turn out to be useful in proving the unmixedness theorem for Cohen--Macaulay rings.
\begin{definition}
    For $I$ an ideal of a Noetherian ring $R$, $I$ is \emph{unmixed} if for any $\mathfrak{p} \in \Ass(R/I)$, $\height(\mathfrak{p}) = \height(I)$.
\end{definition}
\begin{theorem}
    A Noetherian ring $R$ is Cohen--Macaulay if and only if for any ideal $I = (r_1, \cdots, r_k)$ with $\height(I) = k$, $I$ is unmixed.
\end{theorem}

Throughout the paper, we use the symbol \faExternalLink\xspace for external links. Most statements and definitions are accompanied by such a link directly to the source code of the corresponding statement in \mathlib{} or in a branch of our fork \href{https://github.com/Thmoas-Guan/mathlib4_fork/tree/ABS-Criterion-Project-new}{\faExternalLink}. To keep the links stable,  those pointing to our fork target a fixed commit of the master branch, while the others point to the documentation of \mathlib{}. In some code excerpts the dot notation is used: for example, for \lean{I} of type \lean{Ideal R}, the call \lean{Ideal.depth I M} can be shortened to \lean{I.depth M}.

\subsection*{Acknowledgements} We thank Professor Liang Xiao for his preliminary review and for organizing the collaborative work that made this project possible. We are grateful to Dr. Jiedong Jiang for his helpful suggestions during the revision process, which guided multiple major restructurings of the paper. We also thank Dr. Huanhuan Yu for producing the blueprint of the formalization project. We are also deeply grateful to all the contributors and reviewers who have been involved in the development and update process of both \Lean and \mathlib{}, whose efforts have made our work possible. Special thanks go to Professor Jo\"el Riou for developing the derived category and the $\Ext$ functor, which enabled the homological methods. This work is supported in part by the National Key R\&D Program of China grant 2024YFA1014000.

\section{Preliminaries}\label{sec:pre}

In this section, we review some preliminary results that will be used in later proofs. 

\subsection{Associated Primes}\label{subsec:ass}

Recall the definition of associated primes.

\begin{definition}
For an $R$-module $M$, $\mathfrak{p}$ is an associated prime of $M$ if there is $x \in M$ such that the radical of its annihilator is $\mathfrak{p}$.
\end{definition}

At the start of this project, we only had the definition \lean{associatedPrimes}\href{https://leanprover-community.github.io/mathlib4_docs/Mathlib/RingTheory/Ideal/AssociatedPrime/Basic.html#associatedPrimes}{\faExternalLink} and basic facts about associated primes in \mathlib{}. Thanks to the work of Jinzhao Pan, we now have the finiteness of associated primes for finitely generated modules over Noetherian rings in \mathlib{}\href{https://leanprover-community.github.io/mathlib4_docs/Mathlib/RingTheory/Ideal/AssociatedPrime/Finiteness.html#associatedPrimes.finite}{\faExternalLink}. Combined with the fact that the complement of the non-zero-divisors of a module equals the union of its associated primes\href{https://leanprover-community.github.io/mathlib4_docs/Mathlib/RingTheory/Regular/IsSMulRegular.html#biUnion_associatedPrimes_eq_compl_regular}{\faExternalLink}, this lemma allows us to apply the prime avoidance lemma when looking for a regular element in an ideal.

In addition to existing results in \mathlib{}, we formalized results on associated primes of localized modules.
\begin{lemma}
    \href{https://leanprover-community.github.io/mathlib4_docs/Mathlib/RingTheory/Ideal/AssociatedPrime/Localization.html#Module.associatedPrimes.preimage_comap_associatedPrimes_eq_associatedPrimes_of_isLocalizedModule}{\faExternalLink}
    For a multiplicative set $S$ of a Noetherian ring $R$ and a finitely generated $R$-module $M$, denote $R'$ and $M'$ as the localizations of $R$ and $M$ at $S$, one has $\Ass(M) \cap \Spec(R') = \Ass(M')$.
\end{lemma}
\begin{lstlisting}
lemma preimage_comap_associatedPrimes_eq_associatedPrimes_of_isLocalizedModule [IsNoetherianRing R] :
    (Ideal.comap (algebraMap R R')) ⁻¹' (associatedPrimes R M) = associatedPrimes R' M'
\end{lstlisting}
We implemented the left side as the primes whose preimages in $R$ are associated primes of $M$.

%Furthermore, applying the above to the localization at $\mathfrak{p}$, where $\mathfrak{p}$ is a minimal prime over $\Ann(M)$, and since associated primes of $M_p$ exist and any such associated prime must be the maximal ideal of $R_p$ (because $p$ is minimal), we obtain the following lemma.\href{https://leanprover-community.github.io/mathlib4_docs/Mathlib/RingTheory/Ideal/AssociatedPrime/Localization.html#Module.associatedPrimes.minimalPrimes_annihilator_subset_associatedPrimes}{\faExternalLink}

\begin{corollary}\label{lemma:minass}
    For $M$ a finitely generated module over a Noetherian ring $R$, minimal primes over $\Ann(M)$ lie in $\Ass(M)$.
\end{corollary}

\subsection{Krull Dimension of Module}\label{subsec:supportDim}

Recall that the Krull dimension of a module is the Krull dimension of its support viewed as a preordered set. We implement it in \Leanf{} in the same way.\href{https://leanprover-community.github.io/mathlib4_docs/Mathlib/RingTheory/KrullDimension/Module.html#Module.supportDim}{\faExternalLink}

We prove a lemma that will be frequently used later.
\begin{lemma}\label{lemma:quotdim}
    \href{https://leanprover-community.github.io/mathlib4_docs/Mathlib/RingTheory/KrullDimension/Regular.html#Module.supportDim_quotSMulTop_succ_eq_supportDim_mem_jacobson}{\faExternalLink}
    For a Noetherian ring $R$ and a finitely generated $R$-module $M$ and an $M$-regular element $x$ contained in the Jacobson radical of $\Ann(M)$, we have $\dim(M/xM) + 1 = \dim(M)$.
\end{lemma}

\begin{proof}
Note that any $M$-regular element lies outside every minimal prime over $\Ann(M)$. The following two propositions suffice.
\end{proof}

\begin{proposition}
For a Noetherian ring $R$ and a finitely generated $R$-module $M$, quotienting by an element in the Jacobson radical of $\Ann(M)$ decreases the dimension by at most one.
\end{proposition}
\begin{proof}
The formalization of this proposition follows \cite[Proposition 2.3.4]{haution_homalg} via \cite[Lemma 2.3.3]{haution_homalg}\href{https://leanprover-community.github.io/mathlib4_docs/Mathlib/RingTheory/Spectrum/Prime/LTSeries.html#PrimeSpectrum.exist_ltSeries_mem_one_of_mem_last}{\faExternalLink}.
\end{proof}

This can be viewed as a module version of Krull's height theorem for the dimension of quotients. However, since Krull's height theorem is not directly applicable here, we have to argue by considering chains of primes.

\begin{proposition}
For a finitely generated $R$-module $M$, quotienting by an element not in any minimal prime over $\Ann(M)$ decreases the dimension by at least one.
\end{proposition}

\subsection{Projective Dimension}\label{subsec:projdim}

Since the $\Ext$ functor was introduced to \mathlib{} relatively recently as presented in \cite{derived_cat}, at the time of this work, besides the long exact sequence and basic constructions inherited from shifted hom in the derived category such as composition (cup product), we only had a predicate describing that an object $X$ of an abelian category $C$ has projective dimension strictly less than some $n \in \mathbb{N}$, defined as the vanishing of all $\Ext^i(X,Y)$ for all $Y$ and all $i \geq n$ (formalized by Professor Jo\"el Riou).

With only this definition, some operations are awkward (for example, adding a number to the projective dimension in the Auslander--Buchsbaum formula). Even if we try to pick out the minimal $n$ such that  \lean{HasProjectiveDimensionLT X (n + 1)}, the disadvantages are significant: first, such a number may not exist at all (if $X$ does not have a finite projective resolution); second, finding it with \lean{Nat.find} complicates operations and requires decidable instances. For these reasons, we want the projective dimension as a concrete numeric object.

However, the projective dimension cannot be a plain natural number: first, we need to give objects without finite projective resolution the value $\infty$; second, by common convention the projective dimension of the zero object is $-\infty$. We therefore define the projective dimension in \lean{WithBot ℕ∞}, i.e. natural numbers with infinity and negative infinity:
\begin{lstlisting}
def projectiveDimension {C : Type u} [Category.{v, u} C] [Abelian C] (X : C) : WithBot ℕ∞ :=
  sInf {n : WithBot ℕ∞ | ∀ (i : ℕ), n < i → HasProjectiveDimensionLT X i}
\end{lstlisting}
Although this definition looks somewhat verbose, it behaves according to the usual conventions. (Outside code excerpts we denote \lean{projectiveDimension} by $\projdim$.)

With this definition, we can prove the equivalence \lean{projectiveDimension X < n ↔ HasProjectiveDimensionLT X n} and many other useful lemmas to make the numeric version convenient in practice. Many lemmas related to \lean{HasProjectiveDimensionLT} now have corresponding statements for \lean{projectiveDimension}.

We also prove a lemma specific to the category of $R$-modules that is useful when reducing global dimension to the local case.
\begin{lemma}\label{lemma:projdimeqsup}
    \href{https://leanprover-community.github.io/mathlib4_docs/Mathlib/RingTheory/LocalProperties/ProjectiveDimension.html#ModuleCat.projectiveDimension_eq_iSup_localizedModule_prime}{\faExternalLink}
    For a finitely generated module $M$ over a Noetherian ring $R$,
    $$\projdim(M) = \mathop{\mathrm{sup}}\limits_{\mathfrak{p} \in \Spec(R)} \{\projdim_{R_\mathfrak{p}}(M_\mathfrak{p})\}.$$
\end{lemma}
\begin{lstlisting}
lemma projectiveDimension_eq_iSup_localizedModule_prime (M : ModuleCat.{v} R) [Module.Finite R M] : 
    projectiveDimension M = iSup (p : PrimeSpectrum R), projectiveDimension (M.localizedModule p.1.primeCompl)
\end{lstlisting}
Here due to universe issues, the localized module would be in the universe \lean{max u v}, where \lean{u} is the universe level of the ring $R$. Thankfully, if $R$ is assmued to be small, we can shrink it into universe \lean{v} and form an exact functor.

\begin{proof}
For convenience in formalization, we did not directly prove that $\Ext$ commutes with localization. Noting that the left-hand side being $-\infty$ is equivalent to $M = 0$ and similarly for the right-hand side being $-\infty$, we are only left to show $\projdim_R(M) \leq n$ if and only if $\projdim_{R_{\mathfrak{p}}}(M_\mathfrak{p}) \leq n$ for all primes $\mathfrak{p}$.

We prove by induction on $n$. We construct an exact sequence $0 \rightarrow N \rightarrow F \rightarrow M \rightarrow 0$ with $F$ free, localizing at a prime $\mathfrak{p}$ yields an exact sequence $0 \rightarrow N_\mathfrak{p} \rightarrow F_\mathfrak{p} \rightarrow M_\mathfrak{p} \rightarrow 0$. Using the existing lemma about \lean{HasProjectiveDimensionLT} in short exact sequences, we reduce the case $n + 1$ to the case $n$. We are then left to show that $M$ is projective if and only if each localization $M_\mathfrak{p}$ at a prime is projective, which is already in \mathlib{}.
\end{proof}

The trick of doing dimension shifting to two or more modules simultaneously makes us able to reduce the general equation of projective dimensions to exactness and equivalence of being projective. This trick is also carried out in the following lemma and Lemma \ref{lemma:projdimquot}.

\begin{lemma}\label{lemma:projdimle}
    \href{https://leanprover-community.github.io/mathlib4_docs/Mathlib/RingTheory/LocalProperties/ProjectiveDimension.html#ModuleCat.projectiveDimension_le_projectiveDimension_of_isLocalizedModule}{\faExternalLink}
    For a multiplicative set $S$ of $R$ and an $R$-module $M$, $\projdim_{S^{-1}R}(S^{-1}M) \leq \projdim_R(M)$.
\end{lemma}

By applying the same trick again, we can prove that projective dimension is invariant under arbitrary (semi)linear equivalences across different universe levels. As we know \lean{projectiveDimension} is stated using categorical language, we only know it is invariant under linear equivalence within the same universe level in the first place, but it turns out true for general ones.\href{https://leanprover-community.github.io/mathlib4_docs/Mathlib/Algebra/Category/ModuleCat/ProjectiveDimension.html#ModuleCat.projectiveDimension_eq_of_semiLinearEquiv}{\faExternalLink} Consider $M \cong M'$ via a ring isomorphism $R \cong R'$, then $\bigoplus\limits_{m \in M} R \cong \bigoplus\limits_{m' \in M'} R'$ and the two projections $\bigoplus\limits_{m \in M} \rightarrow M$ and $\bigoplus\limits_{m' \in M'} \rightarrow M'$ form a commuting square, thus having isomorphic kernels. By doing this inductively one can reduce to equivalence of projectivity, which is clear.

\subsection{Global Dimension}\label{subsec:gldimpre}

Recall that the global dimension of a ring is the supremum of projective dimensions of all $R$-modules. With the projective dimension as a numeric object, we implement this definition directly.\href{https://github.com/Thmoas-Guan/mathlib4_fork/blob/ABS-Criterion-Project-new/Mathlib/RingTheory/GlobalDimension.lean#L53-L56}{\faExternalLink}
\begin{lstlisting}
noncomputable def globalDimension : WithBot ℕ∞ := iSup (M : ModuleCat.{v} R), projectiveDimension.{v} M
\end{lstlisting}

Outside code excerpts we denote \lean{globalDimension} by $\gldim$. Note that currently the definition of \lean{globalDimension} depends on the universe of the category, but we will prove that it is invariant when the ring is assumed to be small at the end of this subsection.

The supremum over all modules can actually be reduced to the supremum over all finitely generated modules. To prove this we formalized the following form of \cite[Lemma 18.2]{mats_commalg_2ed}.
\begin{lemma}\label{lemma:gldimletfae}
    \href{https://github.com/Thmoas-Guan/mathlib4_fork/blob/ABS-Criterion-Project-new/Mathlib/RingTheory/GlobalDimension.lean#L69-L91}{\faExternalLink} 
    For a commutative ring $R$ and a non-negative integer $n$, the following are equivalent:
    \begin{enumerate}
        \item $\gldim(R) \leq n$,
        \item $\projdim(M) \leq n$ for all finitely generated $R$-modules $M$,
        \item $\Ext_R^i(M, N)$ vanishes for all $R$-modules $M, N$ for $i > n$,
        \item $\injdim(M) \leq n$ for all $R$-modules $M$.
    \end{enumerate}
\end{lemma}

\begin{proof}
$(1) \rightarrow (2)$, $(3) \leftrightarrow (4)$ and $(4) \rightarrow (1)$ are trivial by definition. For $(2) \rightarrow (3)$, the following proposition suffices.
\end{proof}
\begin{proposition}
    \href{https://github.com/Thmoas-Guan/mathlib4_fork/blob/ABS-Criterion-Project-new/Mathlib/Algebra/Category/ModuleCat/Baer.lean#L93-L113}{\faExternalLink}
    For an $R$-module $M$ and a non-negative integer $n$, if $\Ext_R^{n + 1}(R/I, M)$ vanishes for all ideals $I$ of $R$, then $\injdim(M) \leq n$, i.e. for any $R$-module $N$, $\Ext_R^{n + 1}(N, M)$ vanishes.
\end{proposition}

\begin{proof}
We prove by induction on $n$: using dimension-shifting on $M$, we reduce to the case $n + 1$ to that of $n$. For the case $n = 0$, it is just a reformulation of the existing Baer criterion for injective modules.
\end{proof}

Using Lemma ~\ref{lemma:gldimletfae}, we can obtain the following result as a corollary.
\begin{lemma}[{\cite[Theorem 19.1]{Eisenbud2022}}]\label{lemma:gldimfg}
    \href{https://github.com/Thmoas-Guan/mathlib4_fork/blob/ABS-Criterion-Project-new/Mathlib/RingTheory/GlobalDimension.lean#L93-L109}{\faExternalLink}
    For a Noetherian ring $R$, $\gldim(R) = \mathop{\mathrm{sup}}\limits_{M \ f.g.} \{ \projdim(M) \}$.
\end{lemma}

Combining this with the localization results from the previous subsection, we formalize the following.
\begin{lemma}\label{lemma:gldimeqsup}
    \href{https://github.com/Thmoas-Guan/mathlib4_fork/blob/ABS-Criterion-Project-new/Mathlib/RingTheory/GlobalDimension.lean#L129-L150}{\faExternalLink}
    For a Noetherian ring $R$, $\gldim(R) = \mathop{\mathrm{sup}}\limits_{\mathfrak{p} \in \Spec(R)} \{\gldim (R_\mathfrak{p})\}$.
\end{lemma}

\begin{proof}
Since the left-hand side equals the supremum of projective dimensions of finitely generated modules, by Lemma ~\ref{lemma:projdimeqsup}, we know that the left-hand side is no greater than the right-hand side. Conversely, any $R_\mathfrak{p}$-module $N$ is already an $R$-module with $N = N_\mathfrak{p}$; Lemma ~\ref{lemma:projdimle} yields the reverse inequality, completing the proof.
\end{proof}

Due to the statement of Lemma ~\ref{lemma:projdimeqsup}, currently we can only restrict to \lean{Localization} rather than rings \lean{IsLocalization}. We need \lean{globalDimension} to be invariant under ring isomorphism to fill in the final gap.

From Lemma ~\ref{lemma:gldimfg} we know that it suffices to consider the supremum over finitely generated modules. Notice that when assuming \lean{Small.\{v\} R} and \lean{Small.\{v'\} R}, there is a correspondence between finitely generated modules in \lean{ModuleCat.\{v\} R} and \lean{ModuleCat.\{v'\} R}. The \lean{globalDimension} is indeed invariant when assuming the ring is small.

\section{Definition of Depth}\label{sec:depth}

In this section, we formalize the definition of depth as the minimal index of non-vanishing $\Ext$. First, to relate depth to the maximal length of regular sequences contained in a given ideal, we formalize Rees' theorem (Theorem~\ref{thm:Rees}). We then establish basic properties of depth, such as the inequalities associated with short exact sequences and the fact that quotienting by a regular element decreases depth by exactly one. Finally, we present our formalizations of two classical results about depth: the Auslander--Buchsbaum theorem in Subsection ~\ref{subsec:ABthm} and Ischebeck's theorem in Subsection ~\ref{subsec:Ischebeck}.

Throughout this section, we let $R$ be a Noetherian ring and $M$, $N$ be $R$-modules.

We first formalized Rees' theorem, which relates the vanishing of certain $\Ext$-groups and the existence of regular sequences.

\begin{theorem}[Rees, {\cite[Theorem 28]{mats_commalg_2ed}}]\label{thm:Rees}
\href{https://leanprover-community.github.io/mathlib4_docs/Mathlib/RingTheory/Depth/Rees.html#ModuleCat.exists_isRegular_tfae}{\faExternalLink}
Let $M$ be a finitely generated $R$-module and let $I$ be an ideal of $R$ with $IM \neq M$. Fix an integer $n >0$. The following are equivalent.
\begin{enumerate}
\item $\Ext^i_R(N,M)=0$ for all $i<n$ and every finitely generated $R$-module $N$ with $\Supp(N) \subseteq V(I)$;
\item $\Ext_R^i(R/I, M) =0$ for all $i<n$;
\item there exists a finitely generated $R$-module $N$ with $\Supp(N) = V(I)$ such that $\Ext^i_R(N, M) =0 $ for every $i<n$;
\item there exists an $M$-regular sequence $x_1, \dots, x_n$ of length $n$ in $I$.
\end{enumerate}
\end{theorem}

\subsection{Definition of Depth}\label{subsec:depth}

Now we have Rees' theorem building a bridge between the homological characterization and the existence of regular sequences. We want a definition of $\depth$ that is meaningful both in commutative algebra and homological terms. We therefore adopt the following definition for the minimal order of non-vanishing $\Ext$-groups (often called "grade" in some texts).\href{https://github.com/Thmoas-Guan/mathlib4_fork/blob/ABS-Criterion-Project-new/Mathlib/RingTheory/Regular/Depth.lean#L388-L390}{\faExternalLink}
\begin{lstlisting}
def moduleDepth (N M : ModuleCat.{v} R) : ℕ∞ :=
  sSup {n : ℕ∞ | ∀ i : ℕ, i < n → Subsingleton (Ext N M i)}
\end{lstlisting}
Here $R$ is a commutative ring of \lean{Type u} that is \lean{Small.\{v\}} to ensure \lean{ModuleCat.\{v\} R} has enough projectives. These assumptions on $R$ are used throughout in this section. (Outside code excerpts, \lean{moduleDepth M N} will be denoted by $\depth(M,N)$.)

This definition may look unusual but matches the intended idea: it is the minimal index of non-vanishing $\Ext_R^i(N,M)$, and it equals $\infty$ when no such index exists.

We now specialize to the cases of primary interest. The $\depth$ of a module $M$ with respect to an ideal $I$ is defined by taking the module on the left to be $R/I$, denoted by $\depth_I(M)$, formalized as follows:\href{https://github.com/Thmoas-Guan/mathlib4_fork/blob/ABS-Criterion-Project-new/Mathlib/RingTheory/Regular/Depth.lean#L392-L395}{\faExternalLink}
\begin{lstlisting}
def Ideal.depth (I : Ideal R) (M : ModuleCat.{v} R) : ℕ∞ :=
  moduleDepth (ModuleCat.of R (Shrink.{v} (R ⧸ I))) M
\end{lstlisting}
For a local ring $R$, $\depth(M)$ denotes the $\depth$ of $M$ with respect to the unique maximal ideal.\href{https://github.com/Thmoas-Guan/mathlib4_fork/blob/ABS-Criterion-Project-new/Mathlib/RingTheory/Regular/Depth.lean#L397-L399}{\faExternalLink}
\begin{lstlisting}
def IsLocalRing.depth [IsLocalRing R] (M : ModuleCat.{v} R) : ℕ∞ := (IsLocalRing.maximalIdeal R).depth M
\end{lstlisting}

Using Rees' theorem, we relate the above homological notion to the maximal length of $M$-regular sequences contained in a given ideal. In the statements below, we assume $M$ and $N$ are non-zero, finitely generated $R$-modules.

\begin{proposition}\label{lemma:moduleDepth}
    \href{https://github.com/Thmoas-Guan/mathlib4_fork/blob/ABS-Criterion-Project-new/Mathlib/RingTheory/Regular/Depth.lean#L467-L486}{\faExternalLink}
    If $I$ is an ideal of $R$ with $IM < M$ and $\Supp(N) = V(I)$, then $\depth(N,M) = \depth_I(M)$.
\end{proposition}
\begin{proof}
    This is a direct application of the implications $(3) \rightarrow (2) \rightarrow (1)$ of Theorem ~\ref{thm:Rees}.
\end{proof}
This result relates the general definition of $\depth(M,N)$ with the usual $\depth_I(N)$.

With the same hypothesis on $I$, we have the following relation with the maximal length of regular sequences.
\begin{proposition}
    \href{https://github.com/Thmoas-Guan/mathlib4_fork/blob/ABS-Criterion-Project-new/Mathlib/RingTheory/Regular/Depth.lean#L634-L655}{\faExternalLink}
    $\depth(N,M)$ equals the supremum of lengths of $M$-regular sequences contained in $I$.
\end{proposition}
%\begin{lstlisting}
%lemma moduleDepth_eq_sSup_length_regular : moduleDepth N M = sSup {(List.length rs : ℕ∞) |
%    (rs : List R) (_ : RingTheory.Sequence.IsRegular M rs) (_ : ∀ r ∈ rs, r ∈ I)}
%\end{lstlisting}
\begin{proof}
    This is a direct application of the implications $(3) \rightarrow (4) \rightarrow (1)$ of Theorem ~\ref{thm:Rees}.
\end{proof}

Returning to the homological definition, we obtain the following inequalities for $\depth$ in short exact sequences.
\begin{lemma}[A generalization of {\cite[00LX]{stacks-project} }]\label{lemma:depth-ineq}
    \href{https://github.com/Thmoas-Guan/mathlib4_fork/blob/ABS-Criterion-Project-new/Mathlib/RingTheory/Regular/Depth.lean#L582-L632}{\faExternalLink}
    For a short exact sequence of non-zero finitely generated $R$-modules $0 \rightarrow M' \rightarrow M \rightarrow M'' \rightarrow 0$, we have
    \begin{enumerate}
        \item $\depth(N,M) \geq min \{ \depth(N,M'), \depth(N,M'') \}$,
        \item $ \depth(N,M'') \geq min \{ \depth(N,M), \depth(N,M') - 1 \}$,
        \item $ \depth(N,M') \geq min \{ \depth(N,M), \depth(N,M'') + 1 \}$.
    \end{enumerate}
\end{lemma}
\begin{proof}
These follow readily from the long exact sequence.
\end{proof}

Indeed, we mainly use the version of depth with respect to an ideal, for example in Subsection \ref{subsec:gldim} when considering the $\depth$ of a syzygy in a free resolution. However, we will see how the general $\depth$ between two modules and the general form of this inequation are useful in Ischbeck's Theorem (Theorem \ref{thm:Ischbeck}).

We also formalized the behavior of depth under quotient by a regular element.
\begin{theorem}[A generalization of {\cite[090R]{stacks-project}}]\label{thm:quotdepth}
    \href{https://github.com/Thmoas-Guan/mathlib4_fork/blob/ABS-Criterion-Project-new/Mathlib/RingTheory/Regular/Depth.lean#L742-L790}{\faExternalLink}
    For an $M$-regular element $x \in \Ann(N)$, we have
    $$\depth(N,M/xM) + 1 = \depth(N,M).$$
\end{theorem}
Note that this is not immediate from the regular sequence characterization of $\depth$.
\begin{proof}
Consider the long exact sequence induced by the short exact sequence $0 \rightarrow M \xrightarrow{\cdot x} M \rightarrow M/xM \rightarrow 0$. Since $x \in Ann(N)$, the maps $\Ext_R^i(N, M) \xrightarrow{\cdot x} \Ext_R^i(N, M)$ are zero for all $i$. The observation that $\Ext_R^i(N, M/xM)$ vanishes if and only if both $\Ext_R^i(N, M)$ and $\Ext_R^{i + 1}(N, M)$ vanish finishes the proof.
\end{proof}

\subsection{The Auslander--Buchsbaum Theorem}\label{subsec:ABthm}

In this section, we prove the following theorem first introduced in \cite{AuslanderBuchsbaum1957}.

\begin{theorem}[Auslander--Buchsbaum]\label{thm:ABthm}
    \href{https://github.com/Thmoas-Guan/mathlib4_fork/blob/ABS-Criterion-Project-new/Mathlib/RingTheory/Regular/AuslanderBuchsbaum.lean#L314-L436}{\faExternalLink}
    Let $(R, \mathfrak{m})$ be a Noetherian local ring and let $M$ be a non-zero finitely generated $R$-module. If $\projdim(M)$ is finite, then
    $$\projdim(M) + \depth(M) = \depth(R).$$
\end{theorem}
\begin{proof}

The case $\projdim(M)=0$ is immediate since then $M$ is projective and hence free, $\depth(M) = \depth(R)$.

The case $\projdim(M)=1$ is the core, so we prove it separately, following the proof in \cite[Theorem 19.1]{mats_commring}. Pick a $k$-basis of $M/\mathfrak{m}M$ and lift it back to $M$, obtaining a map $f : R^n \rightarrow M$ whose induced map $(R/\mathfrak{m})^n \rightarrow M/\mathfrak{m}M$ is an isomorphism; hence $\Ker(f) \subseteq \mathfrak{m}R^n$. Since $\Ker(f)$ has projective dimension $0$, it is projective, hence free. 

In ordinary proofs, we could view the inclusion $i : \Ker(f) \hookrightarrow R^n$ as a matrix with all elements contained in $\mathfrak{m}$. However, during formalization, viewing a map between free modules as a matrix and operating on it would be complicated. In this case if we follow this setup, we need to take the basis and decompose it into a sum of linear maps corresponding to single matrix elements. So with this observation, instead we set up the following two propositions as stepping stones to help simplify the formalization.

\begin{proposition}\label{lem:vanish1}
    \href{https://github.com/Thmoas-Guan/mathlib4_fork/blob/ABS-Criterion-Project-new/Mathlib/RingTheory/Regular/AuslanderBuchsbaum.lean#L44-L52}{\faExternalLink}
    For a finitely generated free module $M$, if the range of $f : M \rightarrow N$ is contained in $IN$ where $I$ is an ideal of $R$, then $f \in I \cdot (M \rightarrow N)$.
\end{proposition}

\begin{proposition}\label{lem:vanish2}
    \href{https://github.com/Thmoas-Guan/mathlib4_fork/blob/ABS-Criterion-Project-new/Mathlib/RingTheory/Regular/AuslanderBuchsbaum.lean#L171-L186}{\faExternalLink}
    For $R$-modules $L, M, N$, if $f : M \rightarrow N$ is contained in $\Ann(L) \cdot (M \rightarrow N)$, then the map $\Ext_R^i(L,M) \rightarrow \Ext_R^i(L, N)$ induced by $f$ is zero.
\end{proposition}
%\begin{proof}
%    One verifies this lemma easily by reducing to the case $f = a \cdot g$ with $a \in \Ann(L)$.
%\end{proof}

Then we know $i : \Ker(f) \hookrightarrow R^n$ is contained in $\mathfrak{m} \cdot (\Ker(f) \rightarrow R^n)$ from Proposition \ref{lem:vanish1}, so Proposition \ref{lem:vanish2} shows that the map $\Ext_R^i(k, \Ker(f)) \rightarrow \Ext_R^i(k, R^n)$ induced by $i : \Ker(f) \hookrightarrow R^n$ is zero.

A direct calculation on the long exact sequence induced by $0 \rightarrow \Ker(f) \rightarrow R^n \rightarrow M \rightarrow 0$ gives $1 + \depth(M) = \depth(R)$.

Finally, we prove the theorem by induction on $\projdim(M)$. Consider the case $\projdim(M)=n + 1$. Take a surjection $f : R^m \rightarrow M$. By the induction hypothesis, $n + \depth(\Ker(f)) = \depth(R)$. What differs from earlier cases is that $\depth(R^m) = \depth(R) \geq n > 0$, there are no non-zero $k \rightarrow \Ker(f) \hookrightarrow R^m$, thus $\Ext_R^0(k, \Ker(f))$ vanishes.

Again, considering the long exact sequence induced by $0 \rightarrow \Ker(f) \rightarrow R^n \rightarrow M \rightarrow 0$ we know that $\depth(\Ker(f)) = \depth(M) + 1$, again using the fact $n > 0$. Combining this with induction hypothesis finishes the proof.
\end{proof}

\subsection{The Ischebeck Theorem}\label{subsec:Ischebeck}

With the general definition \lean{moduleDepth}, we can state Ischebeck's theorem (first introduced in \cite{Ischebeck1969}) as follows.
\begin{theorem}[Ischbeck]\label{thm:Ischbeck}
    \href{https://github.com/Thmoas-Guan/mathlib4_fork/blob/ABS-Criterion-Project-new/Mathlib/RingTheory/Regular/Ischebeck.lean#L36-L236}{\faExternalLink}
    For a Noetherian local ring $R$ and non-zero finitely generated modules $M$ and $N$, $\depth(N,M) \geq \depth(M) - \dim(N)$.
\end{theorem}
\begin{lstlisting}
theorem moduleDepth_ge_depth_sub_dim : moduleDepth N M ≥
    IsLocalRing.depth M - (Module.supportDim R N).unbot (Module.supportDim_ne_bot_of_nontrivial R N)
\end{lstlisting}
Although we usually avoid subtraction and the operation of \lean{WithBot.unbot}, the proof here actually benefited from this design. Leaving a single \lean{moduleDepth} on the left side allows us to directly unfold and simplify it to conditions on the vanishing of certain $\Ext$-groups. 
\begin{proof}
Since the case $\dim(N) = \infty$ is impossible, we prove the statement by strong induction on $\dim(N)$. If $\dim(N) = 0$, then $\Supp(N) = V(\mathfrak{m})$, and hence $\depth(N,M) = \depth(M)$ by Proposition ~\ref{lemma:moduleDepth}.

Now assume the statement holds for modules of dimension no greater than $n$, and let $\dim(N)=n+1$. In the development of associated primes, one has an induction principle for finitely generated modules over a Noetherian ring: if a condition is stable under extension by exact sequence, and true for the zero module and all quotients of $R$ by prime ideals, then it is true for all finitely generated modules.\href{https://leanprover-community.github.io/mathlib4_docs/Mathlib/RingTheory/Ideal/AssociatedPrime/Finiteness.html#IsNoetherianRing.induction_on_isQuotientEquivQuotientPrime}{\faExternalLink} Preservation under extension by short exact sequences follows from the dual version of Lemma ~\ref{lemma:depth-ineq}\href{https://github.com/Thmoas-Guan/mathlib4_fork/blob/ABS-Criterion-Project-new/Mathlib/RingTheory/Regular/Depth.lean#L530-L542}{\faExternalLink} and a straightforward calculation. This is also an advantage of having the general depth between two modules. Care must be taken to apply the induction hypothesis appropriately to $N_1$ in the two subcases when $\dim(N_1) \leq n$ and when $\dim(N_1) = n + 1$, because in this process the cited induction principle only provides hypotheses for modules of dimension $n + 1$. The hypotheses in the former case come from the strong induction on dimension (and similarly for $N_3$). Thus the remaining and essential case is $N\cong R/\mathfrak{p}$, which is handled by the same argument as in \cite[Lemma 15.2]{mats_commalg_2ed}.
\end{proof}

\begin{corollary}[{\cite[Theorem 29]{mats_commalg_2ed}}]\label{lemma:leass}
    \href{https://github.com/Thmoas-Guan/mathlib4_fork/blob/ABS-Criterion-Project-new/Mathlib/RingTheory/Regular/Ischebeck.lean#L243-L262}{\faExternalLink}
    For a finitely generated $R$-module $M$ and any $\mathfrak{p} \in \Ass(M)$, $\depth(M) \leq \dim(R/\mathfrak{p})$.
\end{corollary}
\begin{proof}
Since $\depth(R/\mathfrak{p},M) = 0$, by Theorem ~\ref{thm:Ischbeck}, we have $\depth(M) \leq \dim(R/\mathfrak{p})$.
\end{proof}

Moreover, from this we deduce the well-known inequalities.
\begin{corollary}\label{cor:depthfin}
    \href{https://github.com/Thmoas-Guan/mathlib4_fork/blob/ABS-Criterion-Project-new/Mathlib/RingTheory/Regular/Ischebeck.lean#L264-L284}{\faExternalLink}
    For a finitely generated non-zero $R$-module $M$, we have $\depth(M) \leq \dim(M) \leq \dim(R)$, all being finite.
\end{corollary}
Although this can be directly proven using Krull's height theorem and the regular sequence characterization, the stronger result above is important in the theory of Cohen--Macaulay modules.

\section{Cohen--Macaulay Modules and Rings}\label{sec:CMdef}

In this section, we first define Cohen--Macaulay modules over local rings. We then prove that passing to the quotient by a regular element preserves the Cohen--Macaulay property, and that localization of a Cohen--Macaulay module remains Cohen--Macaulay. Similarly, we establish corresponding results for Cohen--Macaulay rings. In Subsection ~\ref{subsec:unmixed}, we prove the unmixedness theorem for Cohen--Macaulay rings; as a byproduct, we obtain that a Cohen--Macaulay local ring is catenary. Since a polynomial ring over a Cohen--Macaulay ring is again Cohen--Macaulay (proved in Subsection ~\ref{subsec:ucatenary}), it will be straightforward to deduce that a Cohen--Macaulay ring is universally catenary once the preliminary constructions are in place. In this section, we assume that $R$ is a Noetherian ring.

With a well-developed definition of depth, we define a Cohen--Macaulay module (over a local ring) following \cite[16.A]{mats_commalg_2ed}.
\begin{definition}
    \href{https://github.com/Thmoas-Guan/mathlib4_fork/blob/ABS-Criterion-Project-new/Mathlib/RingTheory/CohenMacaulay/Basic.lean#L51-L54}{\faExternalLink}
    An $R$-module $M$ is a \emph{Cohen--Macaulay module} if $M = 0$ or $\depth(M) = \dim(M)$.
\end{definition}
\begin{lstlisting}
class ModuleCat.IsCohenMacaulay [IsLocalRing R] [Small.{v} R] (M : ModuleCat.{v} R) : Prop where
  depth_eq_dim : Subsingleton M ∨ Module.supportDim R M = IsLocalRing.depth M
\end{lstlisting}

Corollary ~\ref{lemma:leass} yields the following useful lemma, formalized as two separate statements. 
\begin{lemma}[{\cite[Theorem 30(i)]{mats_commalg_2ed}}]\label{lemma:assmin}
    \href{https://github.com/Thmoas-Guan/mathlib4_fork/blob/ABS-Criterion-Project-new/Mathlib/RingTheory/CohenMacaulay/Basic.lean#L86-L96}{\faExternalLink}\href{https://github.com/Thmoas-Guan/mathlib4_fork/blob/ABS-Criterion-Project-new/Mathlib/RingTheory/CohenMacaulay/Basic.lean#L98-L141}{\faExternalLink}
    For a finitely generated Cohen--Macaulay module $M$ over $R$, for any $\mathfrak{p} \in \Ass(M)$, $\dim(M) = \depth(M) = \dim(R/\mathfrak{p})$. Consequently $\mathfrak{p}$ is a minimal prime over $\Ann(M)$.
\end{lemma}
\begin{proof}
    Since $\Ann(M) \subseteq \mathfrak{p}$, by Corollary ~\ref{lemma:leass}, we have $\dim(M) = \depth(M) \leq \dim(R/\mathfrak{p}) \leq \dim(R/\Ann(M)) = \dim(M)$, so every inequality in the chain must be an equality. Hence $\dim(M) = \depth(M) = \dim(R/\mathfrak{p})$, forcing $\mathfrak{p}$ to be a minimal prime over $\Ann(M)$.
\end{proof}

Let $(R, \mathfrak{m})$ be a Noetherian local ring, for a finitely generated $R$-module $M$, quotienting by an $M$-regular element in $\mathfrak{m}$ reduces both the depth and the Krull dimension of $M$ by exactly one; therefore, whether $M$ is Cohen--Macaulay is preserved under such a quotient.
\begin{lemma}\label{lemma:CMquot}
    \href{https://github.com/Thmoas-Guan/mathlib4_fork/blob/ABS-Criterion-Project-new/Mathlib/RingTheory/CohenMacaulay/Basic.lean#L149-L162}{\faExternalLink}
    For a finitely generated module $M$ over $R$ and $x \in \mathfrak{m}$ an $M$-regular element, $M$ is Cohen--Macaulay if and only if $M/xM$ is Cohen--Macaulay.
\end{lemma}
\begin{proof}
    This is a direct consequence of Theorem ~\ref{thm:quotdepth} and Lemma ~\ref{lemma:quotdim}.
\end{proof}

%A similar result holds for quotienting by a regular sequence, which would be used later in Subsection ~\ref{subsec:unmixed}.\href{https://github.com/Thmoas-Guan/mathlib4_fork/blob/ABS-Criterion-Project-new/Mathlib/RingTheory/CohenMacaulay/Basic.lean#L164-L176}{\faExternalLink}
%\begin{lemma}[{\cite[Theorem 30(ii)]{mats_commalg_2ed}}]\label{lemma:CMquotseq}
%	For $M$ a finitely generated  module over $R$ and $(r_1, \cdots, r_k)$ an $M$-regular sequence, then $M$ is Cohen--Macaulay if and only if $M/(r_1, \cdots, r_k)M$ is Cohen--Macaulay.
%\end{lemma}
%\begin{proof}
%	This can be proved easily by induction on the length of the sequence.
%\end{proof}

We then formalize that localization of a Cohen--Macaulay module at a prime ideal is again Cohen--Macaulay.
\begin{lemma}[{\cite[Theorem 30(iii)]{mats_commalg_2ed}}]
    \href{https://github.com/Thmoas-Guan/mathlib4_fork/blob/ABS-Criterion-Project-new/Mathlib/RingTheory/CohenMacaulay/Basic.lean#L378-L388}{\faExternalLink}\href{https://github.com/Thmoas-Guan/mathlib4_fork/blob/ABS-Criterion-Project-new/Mathlib/RingTheory/CohenMacaulay/Basic.lean#L282-L376}{\faExternalLink}
    For a finitely generated Cohen--Macaulay module $M$ over $R$ and a prime ideal $\mathfrak{p}$ of $R$, $M_\mathfrak{p}$ is Cohen--Macaulay over $R_\mathfrak{p}$. Furthermore, if $M_\mathfrak{p} \neq 0$ we have $\depth_\mathfrak{p}(M) = \depth_{R_\mathfrak{p}}(M_\mathfrak{p})$.
\end{lemma}

\begin{proof}
From the regular sequence characterization, we have $\depth_\mathfrak{p}(M) \le \depth(M_\mathfrak{p})$, then Corollary ~\ref{cor:depthfin} implies they are both finite. We proceed by induction on $\depth_\mathfrak{p}(M)$ following the proof in \cite[Theorem 30(iii)]{mats_commalg_2ed}.
\end{proof}

The formalization of a Cohen--Macaulay local ring is done more directly rather than by extending the module notion.
\begin{definition}\label{def:CMlocal}
    \href{https://github.com/Thmoas-Guan/mathlib4_fork/blob/ABS-Criterion-Project-new/Mathlib/RingTheory/CohenMacaulay/Basic.lean#L418-L420}{\faExternalLink}
    A local ring $R$ is a \emph{Cohen--Macaulay local ring} if $\dim(R) = \depth(R)$.
\end{definition}
\begin{lstlisting}
class IsCohenMacaulayLocalRing : Prop extends IsLocalRing R where
  depth_eq_dim : ringKrullDim R = IsLocalRing.depth (ModuleCat.of R R)
\end{lstlisting}
This is equivalent to \lean{(ModuleCat.of R R).IsCohenMacaulay}, but avoids passing through the dimension of modules.\href{https://github.com/Thmoas-Guan/mathlib4_fork/blob/ABS-Criterion-Project-new/Mathlib/RingTheory/CohenMacaulay/Basic.lean#L439-L442}{\faExternalLink}

From the corresponding results on Cohen--Macaulay modules, we obtain the corresponding version for Cohen--Macaulay local rings as a corollary.
\begin{corollary}
    \href{https://github.com/Thmoas-Guan/mathlib4_fork/blob/ABS-Criterion-Project-new/Mathlib/RingTheory/CohenMacaulay/Basic.lean#L447-L455}{\faExternalLink}
    For a Cohen--Macaulay local ring $R$ and $\mathfrak{p}$ a prime ideal of $R$, then $R_\mathfrak{p}$ is Cohen--Macaulay.
\end{corollary}

\begin{corollary}
    \href{https://github.com/Thmoas-Guan/mathlib4_fork/blob/ABS-Criterion-Project-new/Mathlib/RingTheory/CohenMacaulay/Basic.lean#L552-L564}{\faExternalLink}
    For a local ring $(R, \mathfrak{m})$ and an $R$-regular element $x \in \mathfrak{m}$, $R$ is Cohen--Macaulay if and only if $R/(x)$ is Cohen--Macaulay.
\end{corollary}

Due to the change of rings, these are indeed still a small step away from the Cohen--Macaulay module arguments.

With Definition ~\ref{def:CMlocal}, we can define general (non-local) Cohen--Macaulay rings.
\begin{definition}
    \href{https://github.com/Thmoas-Guan/mathlib4_fork/blob/ABS-Criterion-Project-new/Mathlib/RingTheory/CohenMacaulay/Basic.lean#L462-L465}{\faExternalLink}
    A ring $R$ is a \emph{Cohen--Macaulay ring} if $R_\mathfrak{p}$ is a Cohen--Macaulay local ring for every prime ideal $\mathfrak{p}$.
\end{definition}
\begin{lstlisting}
class IsCohenMacaulayRing : Prop where
  CM_localize : ∀ p : Ideal R, ∀ (_ : p.IsPrime), IsCohenMacaulayLocalRing (Localization.AtPrime p)
\end{lstlisting}

Using the above localization result for Cohen--Macaulay local rings, we find that the Cohen--Macaulay property need only be verified at maximal ideals.
\begin{lemma}
    \href{https://github.com/Thmoas-Guan/mathlib4_fork/blob/ABS-Criterion-Project-new/Mathlib/RingTheory/CohenMacaulay/Basic.lean#L475-L502}{\faExternalLink}
    A Noetherian ring $R$ is Cohen--Macaulay if and only if $R_\mathfrak{m}$ is a Cohen--Macaulay local ring for every maximal ideal $\mathfrak{m}$.
\end{lemma}

As a corollary, a Cohen--Macaulay local ring is a Cohen--Macaulay ring\href{https://github.com/Thmoas-Guan/mathlib4_fork/blob/ABS-Criterion-Project-new/Mathlib/RingTheory/CohenMacaulay/Basic.lean#L518-L523}{\faExternalLink} and a local ring that is a Cohen--Macaulay ring is a Cohen--Macaulay local ring \href{https://github.com/Thmoas-Guan/mathlib4_fork/blob/ABS-Criterion-Project-new/Mathlib/RingTheory/CohenMacaulay/Basic.lean#L525-L531}{\faExternalLink}.

\subsection{Unmixedness Theorem for Cohen--Macaulay Rings}\label{subsec:unmixed}

The unmixedness theorem for Cohen--Macaulay rings was first stated in the graded case by Macaulay in \cite{Macaulay1916}. Here we follow the more modern formulation originating from \cite{NorthcottRees1954}. We first prove the following theorem.
\begin{theorem}[{\cite[Theorem 31(iii)]{mats_commalg_2ed}}]\label{thm:htreg}
Let $R$ be a Cohen--Macaulay local ring and let $a_1,\dots,a_r$ be a sequence in the unique maximal ideal $\mathfrak{m}$. The following are equivalent.
\begin{enumerate}
    \item $a_1, \cdots ,a_r$ is $R$-regular,
    \item $\height(a_1, \cdots, a_i) = i$,
        for all $1 \leq i \leq r$,
    \item $\height(a_1, \cdots, a_r) = r$,
    \item $a_1, \cdots ,a_r$ can be extended to a system of parameters.
\end{enumerate}
A system of parameters here means $a_1, \cdots, a_d$ such that $\mathfrak{m}$ is a minimal prime over the ideal they span, where $d = \dim(R)$. 
\end{theorem}

We do not state the theorem as a single \lean{List.TFAE}, but instead provide step-by-step lemmas, since we mainly use the equivalence $(1) \leftrightarrow (3)$\href{https://github.com/Thmoas-Guan/mathlib4_fork/blob/ABS-Criterion-Project-new/Mathlib/RingTheory/CohenMacaulay/Catenary.lean#L42-L74}{\faExternalLink}\href{https://github.com/Thmoas-Guan/mathlib4_fork/blob/ABS-Criterion-Project-new/Mathlib/RingTheory/CohenMacaulay/Catenary.lean#L279-L290}{\faExternalLink} and that some parts do not essentially use the Cohen--Macaulay condition or even the local ones. $(1) \rightarrow (2)$ and $(3) \rightarrow (4)$\href{https://github.com/Thmoas-Guan/mathlib4_fork/blob/ABS-Criterion-Project-new/Mathlib/RingTheory/CohenMacaulay/Catenary.lean#L94-L174}{\faExternalLink} are direct applications of Krull's height theorem. The implication $(4) \rightarrow (1)$ is the only one essentially using the Cohen--Macaulay condition. Then we are able to reduce to the following property about system of parameters.

\begin{proposition}\label{prop:htreg41}
    \href{https://github.com/Thmoas-Guan/mathlib4_fork/blob/ABS-Criterion-Project-new/Mathlib/RingTheory/CohenMacaulay/Catenary.lean#L188-L277}{\faExternalLink}
    For a Cohen--Macaulay local ring $R$, a system of parameters $(r_1, \cdots, r_k)$ forms a regular sequence.
\end{proposition}
\begin{proof}
If the first $k$ elements $(r_1, \cdots, r_k)$ form a regular sequence, then %by Corollary ~\ref{cor:CMquotseq}, 
$R/(r_1, \cdots, r_k)$ is Cohen--Macaulay. As all associated primes of $R/(r_1, \cdots, r_k)$ are minimal by Lemma ~\ref{lemma:assmin}, the fact that $r_{k + 1}$ is not contained in any minimal prime of $R/(r_1, \cdots, r_k)$ implies it is a regular element.
\end{proof}

A by-product of Theorem ~\ref{thm:htreg} is the following lemma.
\begin{lemma}[{\cite[Theorem 31(i), Part1]{mats_commalg_2ed}}]
    \href{https://github.com/Thmoas-Guan/mathlib4_fork/blob/ABS-Criterion-Project-new/Mathlib/RingTheory/CohenMacaulay/Catenary.lean#L321-L327}{\faExternalLink}
    For a Cohen--Macaulay local ring $R$ and $I$ a proper ideal of it, then $\depth_I(R) = \height(I)$.
\end{lemma}
\begin{proof}
By Theorem ~\ref{thm:htreg}, the left-hand side is clearly not greater than the right-hand side. Then as we can find $J \subseteq I$ such that $\height(J) = \height(I)$ and $J$ can be generated by $\height(I)$ elements, Theorem ~\ref{thm:htreg} shows the generators form a regular sequence, and the right-hand side is not greater than the left-hand side.
\end{proof}

We then establish the following theorem which implies that a Cohen--Macaulay local ring is catenary. 
\begin{theorem}[{\cite[Theorem 31(i), Part2]{mats_commalg_2ed}}]
    \href{https://github.com/Thmoas-Guan/mathlib4_fork/blob/ABS-Criterion-Project-new/Mathlib/RingTheory/CohenMacaulay/Catenary.lean#L377-L394}{\faExternalLink}
    For $R$ a Cohen--Macaulay local ring and a proper ideal $I$ of it, we have $\height(I) + \dim(R/I) = \dim(R)$.
\end{theorem}
\begin{proof}
The formalization then follows the proof in \cite[Theorem 31(i), Part 2]{mats_commalg_2ed}. However, to reduce to the case of prime ideals, this should be direct by operating with supremum and infimum, but the implementation met some problems. We can see we are operating with natural numbers with both infinity and negative infinity, this makes many adding rules unusable.

It turns out we are already in a Noetherian ring, minimal primes over $I$ are finite, height as an infimum among them and the dimension of the quotient as a supremum among them can both be attained by a certain prime, so we are done.
\end{proof}

From this result and the fact that the localization of a Cohen--Macaulay ring is again Cohen--Macaulay, for prime ideals $\mathfrak{q} \subseteq \mathfrak{p}$, we have $\dim(R_\mathfrak{p}) = \height(\mathfrak{q}R_\mathfrak{p}) + \dim(R_\mathfrak{p}/\mathfrak{q}R_\mathfrak{p})$, i.e. $\height(\mathfrak{p}) - \height(\mathfrak{q}) = \height(\mathfrak{p}/\mathfrak{q})$. This implies that a Cohen--Macaulay local ring is catenary. However, since the notions of catenary order and catenary rings are not yet in \mathlib{}, we do not develop this further here.

Returning to the main theorem, we first define the notion of an unmixed ideal.
\begin{definition}
    \href{https://github.com/Thmoas-Guan/mathlib4_fork/blob/ABS-Criterion-Project-new/Mathlib/RingTheory/CohenMacaulay/Unmixed.lean#L25-L27}{\faExternalLink}
    For an ideal $I$ of a Noetherian ring $R$, $I$ is \emph{unmixed} if for any $\mathfrak{p} \in \Ass(R/I)$, $\height(\mathfrak{p}) = \height(I)$.
\end{definition}
\begin{lstlisting}
class Ideal.IsUnmixed (I : Ideal R) : Prop where
  height_eq : ∀ {p : Ideal R}, p ∈ associatedPrimes R (R ⧸ I) → p.height = I.height
\end{lstlisting}
In particular, if $I = (r_1, \cdots, r_k)$ with $\height(I) = k$, then $I$ is unmixed if and only if $R/I$ has no embedded primes.

\begin{theorem}[Unmixedness Theorem]
    A Noetherian ring $R$ is Cohen--Macaulay if and only if for any ideal $I = (r_1, \cdots, r_k)$ with $\height(I) = k$, $I$ is unmixed.
\end{theorem}
\begin{lstlisting}
theorem isCohenMacaulayRing_iff_unmixed : IsCohenMacaulayRing R ↔
    ∀ (l : List R), (Ideal.ofList l).height = l.length → (Ideal.ofList l).IsUnmixed
\end{lstlisting}
\begin{proof}
With the results earlier in this section, the formalization of the unmixedness theorem follows the proof in \cite[Theorem 32]{mats_commalg_2ed}.\href{https://github.com/Thmoas-Guan/mathlib4_fork/blob/ABS-Criterion-Project-new/Mathlib/RingTheory/CohenMacaulay/Unmixed.lean#L128-L180}{\faExternalLink}
\end{proof}

\subsection{Universally Catenary}\label{subsec:ucatenary}

Using the results proven in Subsection ~\ref{subsec:unmixed}, to prove that Cohen--Macaulay rings are universally catenary, it remains only to show that a polynomial ring over a Cohen--Macaulay ring is also Cohen--Macaulay.\href{https://github.com/Thmoas-Guan/mathlib4_fork/blob/ABS-Criterion-Project-new/Mathlib/RingTheory/CohenMacaulay/Polynomial.lean#L127-L199}{\faExternalLink}

\begin{lemma}
The polynomial ring $R[X]$ over a Cohen--Macaulay ring $R$ is Cohen--Macaulay.
\end{lemma}

Assume $R$ is Cohen--Macaulay. For any prime ideal $\mathfrak{p}$ of $R[X]$, set $\mathfrak{q} = \mathfrak{p}\cap R$. Note that $R[X]_\mathfrak{p}$ is the localization of $R_\mathfrak{q}[X]$ at $\mathfrak{p}R_\mathfrak{q}[X]$, so the problem reduces to the following case: if $R$ is a Cohen--Macaulay local ring and $\mathfrak{p}$ is a prime ideal of $R[X]$ with $\mathfrak{p}\cap R$ equal to the maximal ideal of $R$, then $R[X]_\mathfrak{p}$ is a Cohen--Macaulay local ring.

Although the notions of catenary order and catenary rings are not yet available in \mathlib{}, once they are added, it will be straightforward to combine the above results to conclude that Cohen--Macaulay rings are universally catenary.

\section{Regular Local Ring}\label{sec:regloc}

In this section, we first introduce our formalization of the definition of regular local rings and prove that they are domains. Then we prove the two directions of the Auslander--Buchsbaum--Serre criterion in the next two subsections. We first prove that the global dimension of a regular local ring equals its Krull dimension by showing that a maximal Cohen--Macaulay module over a regular local ring is free. Then we prove the converse: if the unique maximal ideal of a Noetherian local ring has finite projective dimension, then it is generated by a regular sequence; this completes the converse direction and can be viewed as a weakened form of the Ferrand--Vasconcelos theorem. Finally, we establish some corollaries of the Auslander--Buchsbaum--Serre criterion in the last subsection, including that regularity of a ring can be checked at maximal ideals only, and Hilbert's Syzygy Theorem.

We define a regular local ring simply as a Noetherian local ring whose maximal ideal has span rank (minimal number of generators) equal to its Krull dimension.\href{https://github.com/Thmoas-Guan/mathlib4_fork/blob/ABS-Criterion-Project-new/Mathlib/RingTheory/RegularLocalRing/Defs.lean#L36-L39}{\faExternalLink}
\begin{lstlisting}
class IsRegularLocalRing : Prop extends IsLocalRing R, IsNoetherianRing R where
  reg : (maximalIdeal R).spanFinrank = ringKrullDim R
\end{lstlisting}

Using the fact that the span rank of the maximal ideal equals the dimension of the cotangent space over the residue field\href{https://github.com/Thmoas-Guan/mathlib4_fork/blob/ABS-Criterion-Project-new/Mathlib/RingTheory/RegularLocalRing/Defs.lean#L77-L117}{\faExternalLink}, we establish the equivalence of this definition with the cotangent-space formulation.\href{https://github.com/Thmoas-Guan/mathlib4_fork/blob/ABS-Criterion-Project-new/Mathlib/RingTheory/RegularLocalRing/Defs.lean#L119-L121}{\faExternalLink}

\subsection{Regular Local Ring is Domain}\label{subsec:regdom}

To prove that a regular local ring is a domain and that minimal generators of the maximal ideal form a regular sequence, we establish the following lemma.\href{https://github.com/Thmoas-Guan/mathlib4_fork/blob/ABS-Criterion-Project-new/Mathlib/RingTheory/RegularLocalRing/Basic.lean#L50-L195}{\faExternalLink}
\begin{lemma}\label{lemma:regiff}
Let $(R,\mathfrak{m},k)$ be a regular local ring of dimension $d$, and let $a_1,\cdots,a_k \in \mathfrak{m}$. Then the following are equivalent.
\begin{enumerate}
  \item $a_1,\cdots,a_k$ can be extended to a regular system of parameters of $R$.
  \item The images $\overline{a_1},\cdots,\overline{a_k}\in\mathfrak{m}/\mathfrak{m}^2$ are $k$-linearly independent.
  \item $R/(a_1,\cdots,a_k)$ is a regular local ring of dimension $d-k$.
\end{enumerate}
A regular system of parameters here simply means a set of minimal generators of $\mathfrak{m}$.
\end{lemma}

This lemma is about studying the cotangent space of the quotient by a list of elements. However during formalization we need to switch between $R$-scalars and $k$-scalars, and this is where most of the inconvenience comes from.

\begin{corollary}\label{coro:regquot}
    \href{https://github.com/Thmoas-Guan/mathlib4_fork/blob/ABS-Criterion-Project-new/Mathlib/RingTheory/RegularLocalRing/Basic.lean#L197-L202}{\faExternalLink}
    For a regular local ring $R$ and $x \in \mathfrak{m} \setminus \mathfrak{m}^2$, $R/(x)$ is a regular local ring of dimension $\dim(R) - 1$.
\end{corollary}

We then prove that a regular local ring $(R,\mathfrak{m},k)$ is a domain by induction on the dimension.
\begin{theorem}
    \href{https://github.com/Thmoas-Guan/mathlib4_fork/blob/ABS-Criterion-Project-new/Mathlib/RingTheory/RegularLocalRing/Basic.lean#L213-L264}{\faExternalLink}
    A regular local ring $R$ is a domain.
\end{theorem}
\begin{proof}
Proof by induction on $\dim(R)$, dimension zero is clear. For positive dimension, $\mathfrak{m}$ is not minimal, hence by prime avoidance (allowing one of them not a prime) we can pick $x \in \mathfrak{m} \setminus \mathfrak{m}^2$ that is not in any minimal prime of $R$. The induction hypothesis, along with Corollary ~\ref{coro:regquot} shows that $R/(x)$ is a domain. From the property of $x$, we know that $(x)$ is a prime ideal that is not minimal. Consider a minimal prime $\mathfrak{p} \subsetneq (x)$, then $x\mathfrak{p} = \mathfrak{p}$. Nakayama's lemma shows $\mathfrak{p} = 0$, finishing the proof.
\end{proof}

%\cite[Proposition 2.2.3]{CM_ring}.

For a regular local ring $R$ of dimension $d$ and a sequence $rs$ of length $d$ generating its maximal ideal, any quotient of $R$ by a subsequence of $rs$ is still regular and hence a domain by the result above; therefore, $rs$ is an $R$-regular sequence.
\begin{corollary}
    \href{https://github.com/Thmoas-Guan/mathlib4_fork/blob/ABS-Criterion-Project-new/Mathlib/RingTheory/RegularLocalRing/Basic.lean#L266-L329}{\faExternalLink}
    For a regular local ring $R$, any regular system of parameters forms a regular sequence.
\end{corollary}
\begin{corollary}
    \href{https://github.com/Thmoas-Guan/mathlib4_fork/blob/ABS-Criterion-Project-new/Mathlib/RingTheory/CohenMacaulay/Maximal.lean#L59-L78}{\faExternalLink}
    A regular local ring is Cohen--Macaulay.
\end{corollary}

\subsection{Global Dimension of Regular Local Ring}\label{subsec:gldim}

In this section, we prove that the global dimension of a regular local ring equals its Krull dimension.
\begin{theorem}\label{thm:reggldim}
    \href{https://github.com/Thmoas-Guan/mathlib4_fork/blob/ABS-Criterion-Project-new/Mathlib/RingTheory/RegularLocalRing/GlobalDimension.lean#L78-L122}{\faExternalLink}
    For a regular local ring $R$, $\gldim(R) = \dim(R)$.
\end{theorem}
\begin{lstlisting}
theorem IsRegularLocalRing.globalDimension_eq_ringKrullDim [Small.{v} R] [IsRegularLocalRing R] :
    globalDimension.{v} R = ringKrullDim R
\end{lstlisting}

There is a well-known proof of this theorem using the following two facts for a Noetherian local ring $(R, \mathfrak{m}, k)$.
\begin{proposition}[{\cite[Theorem 41]{mats_commalg_2ed}}]\label{prop:projres}
$\gldim(R) = \projdim(k)$.
\end{proposition}
\begin{proposition}[{\cite[Lemma 18.6]{mats_commalg_2ed}}]
For any finitely generated $R$-module $M$ and an $M$-regular element $x$, we have $\projdim(M/xM) = \projdim(M) + 1$.
\end{proposition}

Considering $k$ as the quotient of $R$ by a regular sequence completes the proof. However, the $\Tor$ functor and related infrastructure were still incomplete in \mathlib{} when this work was finished, so we could not prove Proposition \ref{prop:projres} with the current setup. Instead, we prove that maximal Cohen--Macaulay modules over a regular local ring are free and use that to proceed.

\subsubsection{\textbf{Maximal Cohen--Macaulay Module}}

We define a maximal Cohen--Macaulay module over $R$ as a module whose depth equals $\dim(R)$.\href{https://github.com/Thmoas-Guan/mathlib4_fork/blob/ABS-Criterion-Project-new/Mathlib/RingTheory/CohenMacaulay/Maximal.lean#L34-L37}{\faExternalLink}
\begin{definition}
    A module $M$ over a local ring $R$ is a \emph{maximal Cohen--Macaulay module} if $\depth(M) = \dim(R)$.
\end{definition}
\begin{lstlisting}
class ModuleCat.IsMaximalCohenMacaulay [IsLocalRing R] [Small.{v} R] (M : ModuleCat.{v} R) : Prop where
    depth_eq_dim : IsLocalRing.depth M = ringKrullDim R
\end{lstlisting}
As we know from Corollary \ref{cor:depthfin} that for a finitely generated $R$-module $M$, $\depth(M) \leq \dim(M) \leq \dim(R)$, maximal Cohen--Macaulay modules can be viewed as modules over $R$ with the maximum possible $\depth$. In particular, for a finitely generated maximal Cohen--Macaulay module $M$, inequalities in Corollary \ref{cor:depthfin} must be equalities since $\depth(M) = \dim(R)$, hence it is Cohen--Macaulay.\href{https://github.com/Thmoas-Guan/mathlib4_fork/blob/ABS-Criterion-Project-new/Mathlib/RingTheory/CohenMacaulay/Maximal.lean#L49-L57}{\faExternalLink}

\begin{theorem}[{\cite[00NT]{stacks-project}}]\label{thm:maxfree}
    \href{https://github.com/Thmoas-Guan/mathlib4_fork/blob/ABS-Criterion-Project-new/Mathlib/RingTheory/CohenMacaulay/Maximal.lean#L146-L192}{\faExternalLink}
    For a regular local ring $R$, any finitely generated maximal Cohen--Macaulay module over $R$ is free.
\end{theorem}
\begin{proof}
    We proceed by induction on $\dim(R)$. Consider the case $\dim(R)=n + 1$. Choose $x \in \mathfrak{m} \setminus \mathfrak{m}^2$. Since for any $\mathfrak{p}\in\Ass(M)$ we have $\dim(R/\mathfrak{p})\ge\depth(M)=\dim(R)$, this implies $\mathfrak{p} = 0$. $M$ is torsion free, hence $x$ is $M$-regular. Theorem ~\ref{thm:quotdepth} shows $M/xM$ is a maximal Cohen--Macaulay module over $R/(x)$, the following lemma finishes the proof.
\end{proof} 

\begin{lemma}[{\cite[00NS]{stacks-project}}]\label{lemma:free}
    \href{https://leanprover-community.github.io/mathlib4_docs/Mathlib/RingTheory/Regular/Free.html#Module.free_quotSMulTop_iff_free}{\faExternalLink}
    Let $(R, \mathfrak{m}, k)$ be a Noetherian local ring and $M$ a finitely generated $R$-module. Let $x \in \mathfrak{m}$ be $M$-regular. If $M/xM$ is free over $R/(x)$, then $M$ is free over $R$.
\end{lemma}

Coincidentally, this lemma is used again in the next subsection.

%\begin{proof}
%Pick the basis of $M/xM$ over $R/(x)$ and lift it back to $M$, from this we can obtain a surjection $f : R^n \rightarrow M$ with $\Ker(f) \subseteq xR^n$. As $x$ is $M$-regular, $\Ker(f) = x\Ker(f)$ from $\Ker(f) \subseteq xR^n$. By Nakayama lemma, this implies $f$ is an isomorphism.
%\end{proof}

The linked version may appear somewhat different. This is because, by using the general Nakayama's lemma, we can omit the Noetherian local condition and replace it with the assumptions that $x$ lies in the Jacobson radical of $R$ and that $M$ is finitely presented. 

\subsubsection{\textbf{The Main Proof}}

For a regular local ring $(R, \mathfrak{m}, k)$ of dimension $d$ and any $R$-module $M$, repeated application of Lemma ~\ref{lemma:depth-ineq} shows that the $d$-th syzygy of $M$ in a free resolution is either $0$ or maximal Cohen--Macaulay; in either case, it is projective. (See arguments in \cite[Theorem 2.2.7]{CM_ring} for details.) However, formalizing a distant syzygy and the depth calculations is cumbersome, so we use a trick to simplify the formalization.

For a regular local ring $R$, the Auslander--Buchsbaum formula implies that for a finitely generated $R$-module $M$ of finite projective dimension, the projective dimension does not exceed $\depth(R) = \dim(R)$. Since $\depth(k) = 0$, if $k$ has finite projective dimension, then its projective dimension equals $\depth(R) = \dim(R)$. Thus, the following suffices for Theorem ~\ref{thm:reggldim}.

\begin{lemma}
    \href{https://github.com/Thmoas-Guan/mathlib4_fork/blob/ABS-Criterion-Project-new/Mathlib/RingTheory/RegularLocalRing/GlobalDimension.lean#L70-L76}{\faExternalLink}
    For any finitely generated $R$-module $M$, $\projdim(M)$ is finite.
\end{lemma}

The trick is to set up an auxiliary result that is easier to prove by induction on $i \in \mathbb{N}$, using only single-step dimension shifting.
\begin{proposition}
    \href{https://github.com/Thmoas-Guan/mathlib4_fork/blob/ABS-Criterion-Project-new/Mathlib/RingTheory/RegularLocalRing/GlobalDimension.lean#L26-L68}{\faExternalLink}
    For any $i \in \mathbb{N}$ and finietly generated $R$-module $M$, if $\depth(M) + i \geq \dim(R)$, then $\projdim(M)$ is finite.
\end{proposition}
%\begin{lstlisting}
%lemma finite_projectiveDimension_of_isRegularLocalRing_aux [IsRegularLocalRing R] [Small.{v} R]
%    (M : ModuleCat.{v} R) [Module.Finite R M] (i : ℕ) : 
%    IsLocalRing.depth M + i ≥ ringKrullDim R → ∃ n, HasProjectiveDimensionLE M n
%\end{lstlisting}
\begin{proof}
    For $i=0$, $\depth(M)\ge\dim(R)$ implies equality, so $M$ is maximal Cohen--Macaulay and hence free by Theorem ~\ref{thm:maxfree}; thus $M$ has finite projective dimension. Assume the claim for $i$. For $i+1$, suppose $M \neq 0$, take a short exact sequence $0\to N\to F\to M\to 0$ with $F$ free and finitely generated. Assume $N \neq 0$ or else $M$ is free directly, by Lemma ~\ref{lemma:depth-ineq}, $\depth(N) + i \geq \min\{\depth(F), \depth(M) + 1\} + i = \min\{\depth(F) + i, \depth(M) + i + 1\}$. Since $\depth(F)=\depth(R)=\dim(R)$, the right-hand side is at least $\dim(R)$, so $\depth(N)+i \geq \dim(R)$. By the induction hypothesis, $N$ has finite projective dimension; hence, so does $M$. This completes the induction.
\end{proof}

The proof only involves dimension shifting of one term, making it easier to formalize. Combining all the results above yields the formal proof of Theorem ~\ref{thm:reggldim}.

\subsection{A Weakened Version of Ferrand--Vasconcelos Theorem}\label{subsec:FVthm}

The following theorem will directly imply the desired result.
\begin{theorem}[Ferrand--Vasconcelos, {\cite[Theorem 2.2.8]{CM_ring}}]
    For an ideal $I$ of a Noetherian local ring $R$, if $I$ has finite projective dimension and $I/I^2$ is free over $R/I$, then $I$ is generated by a regular sequence.
\end{theorem}
If the residue field $k$ of a Noetherian local ring $(R, \mathfrak{m}, k)$ has finite projective dimension, then so does $\mathfrak{m}$; applying the theorem to $\mathfrak{m}$ shows that $\mathfrak{m}$ is generated by a regular sequence. Since any regular sequence contained in $\mathfrak{m}$ has length at most $\depth(R) \leq \dim(R)$, this implies that $R$ is regular.

However, proving the full Ferrand--Vasconcelos theorem requires showing that any ideal $I$ with a finite free resolution contains an R‑regular element, as shown in \cite[Corollary 1.4.7]{CM_ring}, which in turn requires a theory of rank for free resolutions. As the theory of rank in \mathlib{} is not yet fully developed, the complete proof is currently out of reach. The key observation is that the induction in the proof of the Ferrand--Vasconcelos theorem still works if we only consider the maximal ideal. Thus, we restrict our attention to the unique maximal ideal and prove the result for that case; this already yields regularity when the residue field has finite projective dimension. Let $(R, \mathfrak{m})$ be a Noetherian local ring throughout the following part of this subsection.

\begin{theorem}\label{thm:weakFVthm}
    \href{https://github.com/Thmoas-Guan/mathlib4_fork/blob/ABS-Criterion-Project-new/Mathlib/RingTheory/RegularLocalRing/AuslanderBuchsbaumSerre.lean#L661-L664}{\faExternalLink}
    If $\mathfrak{m}$ has finite projective dimension, then $\mathfrak{m}$ is generated by a regular sequence.
\end{theorem}
\begin{lstlisting}
theorem generate_by_regular [IsLocalRing R] [IsNoetherianRing R] [Small.{v} R]
    (h : ∃ n, HasProjectiveDimensionLE (ModuleCat.of R (Shrink.{v} (maximalIdeal R))) n) :
    ∃ rs : List R, IsRegular R rs ∧ Ideal.ofList rs = maximalIdeal R
\end{lstlisting}

\begin{proof}
The existence of a regular element in the maximal ideal, assuming the maximal ideal has finite projective dimension and $R$ is not a field, follows from the Auslander--Buchsbaum formula (Theorem ~\ref{thm:ABthm}). This allows us to get around the complicated arguments about rank in a free resolution.

\begin{lemma}
    \href{https://github.com/Thmoas-Guan/mathlib4_fork/blob/ABS-Criterion-Project-new/Mathlib/RingTheory/RegularLocalRing/AuslanderBuchsbaumSerre.lean#L294-L357}{\faExternalLink}
    If $R$ is not a field and $\mathfrak{m}$ has finite projective dimension, then $\mathfrak{m}$ contains an $R$-regular element.
\end{lemma}
\begin{proof}
    Suppose there are no regular elements in $\mathfrak{m}$. Then $\depth(R) = 0$. If $k$ has finite projective dimension, the Auslander--Buchsbaum formula (Theorem ~\ref{thm:ABthm}) gives $\projdim(k) + \depth(k) = \depth(R) = 0$. Since $\depth(k) = 0$, we get $\projdim(k) = 0$. The residue field $k$ is projective hence free, thus $\mathfrak{m} = 0$, contradicting the assumption that $R$ is not a field.
\end{proof} 

Using prime avoidance applied to $\mathfrak{m}^2$ and $\Ass(R)$, we can further find an $R$-regular element $x \in \mathfrak{m} \setminus \mathfrak{m}^2$.
\begin{corollary}\label{lemma:existreg}
    \href{https://github.com/Thmoas-Guan/mathlib4_fork/blob/ABS-Criterion-Project-new/Mathlib/RingTheory/RegularLocalRing/AuslanderBuchsbaumSerre.lean#L359-L393}{\faExternalLink}
    If $R$ is not a field and $\mathfrak{m}$ has finite projective dimension, then there exists an $R$-regular element in $\mathfrak{m}\setminus\mathfrak{m}^2$.
\end{corollary}
Returning to the main proof, we proceed by induction on the span rank (minimal number of generators) of $\mathfrak{m}$. Consider the case $n + 1$, by Corollary ~\ref{lemma:existreg}, we can choose an $R$-regular element $x \in \mathfrak{m} \setminus \mathfrak{m}^2$. To apply the induction hypothesis to $R/(x)$, it remains to show its maximal ideal $\Tilde{\mathfrak{m}}$ has finite projective dimension over $R/(x)$. Let $\overline{J}$ be a direct sum complement of $\mathrm{span} \{\overline{x}\}$ in $\mathfrak{m}/\mathfrak{m}^2$, and $J$ be the preimage of $\overline{J}$ in $\mathfrak{m}$. We have $J \cap (x) \subseteq x\mathfrak{m}$, then consider the following maps whose composition is the identity:
$$
\Tilde{\mathfrak{m}} \cong (J + (x))/(x) \cong J/(J \cap (x)) \rightarrow \mathfrak{m}/x\mathfrak{m} \rightarrow \Tilde{\mathfrak{m}}.
$$
Therefore $\Tilde{\mathfrak{m}}$ is a direct summand of $\mathfrak{m}/x\mathfrak{m}$ as an $R/(x)$ module. We are only left to show that $\projdim_{R/(x)}(\mathfrak{m}/x\mathfrak{m})$ is finite. Then the following lemma suffices.
\end{proof}

\begin{lemma}[{\cite[Lemma 1.3.5]{CM_ring}}]\label{lemma:projdimquot}
    \href{https://github.com/Thmoas-Guan/mathlib4_fork/blob/ABS-Criterion-Project-new/Mathlib/RingTheory/RegularLocalRing/AuslanderBuchsbaumSerre.lean#L193-L290}{\faExternalLink}
    For a finitely generated module $M$ over a Noetherian local ring $R$ and for $x \in \mathfrak{m}$ that is both $R$-regular and $M$-regular, we have  $\projdim_R(M) = \projdim_{R/(x)}(M/xM)$.
\end{lemma}
\begin{proof}
Since projective and free coincide over a local ring, Lemma~\ref{lemma:free} implies that $M$ is projective over $R$ if and only if $M/xM$ is projective over $R/(x)$. Using dimension shifting, it suffices to reduce the general equality of projective dimensions to the equivalence of projectivity on both sides. Concretely, for non-projective $M$, take a short exact sequence $0 \rightarrow N \xrightarrow{f} R^n \xrightarrow{g} M \rightarrow 0$. Tensoring with $R/(x)$ yields an exact sequence $N/xN \rightarrow (R/(x))^n \rightarrow M/xM \rightarrow 0$ with the obvious transition maps. Any element $y$ in the kernel of the map $N \rightarrow R^n \rightarrow (R/(x))^n$ lifts to an element of $N$ of the form $y = x \cdot z$. Since $x \cdot g(z) = g (x \cdot z) = 0$, we have $z \in N$. Thus the kernel of $N \rightarrow N/xN \rightarrow (R/(x))^n$ equals $xN$, so the induced map $N/xN \rightarrow (R/(x))^n$ is injective. Because $\projdim_R(N) = \projdim_R(M) - 1$ and $\projdim_{R/(x)} \ N/xN = \projdim_{R/(x)} \ M/xM - 1$, and since $x$ remains $R$-regular and $N$-regular (as $N$ is a submodule of a free module), the equality of projective dimensions follows. 
\end{proof}

A remark about the exactness of $0 \rightarrow N/xN \rightarrow (R/(x))^n \rightarrow M/xM \rightarrow 0$ is that it can be implied by $\Tor_R^1(R/(x), M) = 0$, which can be computed using the free resolution $\cdots 0 \rightarrow R \xrightarrow{\cdot x} R \rightarrow R/(x) \rightarrow 0$. But again, since $\Tor$ is incomplete, we compute exactness directly.

We conclude this subsection with the following two direct corollaries of Theorem \ref{thm:weakFVthm}.
\begin{theorem}
    \href{https://github.com/Thmoas-Guan/mathlib4_fork/blob/ABS-Criterion-Project-new/Mathlib/RingTheory/RegularLocalRing/AuslanderBuchsbaumSerre.lean#L666-L685}{\faExternalLink}
    For a Noetherian local ring $(R, \mathfrak{m})$, if $\mathfrak{m}$ has finite projective dimension then $R$ is regular.
\end{theorem}
\begin{theorem}
    \href{https://github.com/Thmoas-Guan/mathlib4_fork/blob/ABS-Criterion-Project-new/Mathlib/RingTheory/RegularLocalRing/AuslanderBuchsbaumSerre.lean#L687-L694}{\faExternalLink}
    If a Noetherian local ring has finite global dimension, then it is regular.
\end{theorem}

\subsection{The Final Conclusions and Related Results}\label{subsec:concl}

\subsubsection{\textbf{Auslander--Buchsbaum--Serre Criterion}}

With the results in Subsections ~\ref{subsec:gldim} and ~\ref{subsec:FVthm}, we can finally state and prove the Auslander--Buchsbaum--Serre criterion, following the original version in \cite{AuslanderBuchsbaum1956}, as follows.
\begin{theorem}
    \href{https://github.com/Thmoas-Guan/mathlib4_fork/blob/ABS-Criterion-Project-new/Mathlib/RingTheory/RegularLocalRing/Localization.lean#L26-L30}{\faExternalLink}
    A Noetherian local ring $R$ is regular if and only if $\gldim(R)$ is finite.
\end{theorem}

\begin{lstlisting}
theorem Auslander_Buchsbaum_Serre [IsLocalRing R] [IsNoetherianRing R] [Small.{v} R] :
    IsRegularLocalRing R ↔ globalDimension.{v} R < Top.top
\end{lstlisting}

Now we can obtain the following non-homological corollary, proven by this homological criterion.
\begin{corollary}
    \href{https://github.com/Thmoas-Guan/mathlib4_fork/blob/ABS-Criterion-Project-new/Mathlib/RingTheory/RegularLocalRing/Localization.lean#L43-L47}{\faExternalLink}
    For a regular local ring $R$ and $\mathfrak{p}$ a prime ideal of $R$, the localization $R_\mathfrak{p}$ is regular.
\end{corollary}
\begin{proof}
From Lemma ~\ref{lemma:gldimeqsup} we know that $\gldim(R_\mathfrak{p}) \leq \gldim(R)$ must be finite.
\end{proof}

\subsubsection{\textbf{Regular Ring}}

We also developed the basic theory of regular rings.
\begin{definition}
    \href{https://github.com/Thmoas-Guan/mathlib4_fork/blob/ABS-Criterion-Project-new/Mathlib/RingTheory/RegularLocalRing/RegularRing/Basic.lean#L25-L28}{\faExternalLink}
    A Noetherian ring $R$ is a \emph{regular ring} if $R_\mathfrak{p}$ is a regular local ring for every prime ideal $\mathfrak{p}$.
\end{definition}
\begin{lstlisting}
class IsRegularRing : Prop where
  localization_isRegular : ∀ p : Ideal R, ∀ (_ : p.IsPrime), IsRegularLocalRing (Localization.AtPrime p)
\end{lstlisting}

We have the analogous statement of Theorem ~\ref{thm:reggldim} (global dimension argument) for regular rings.
\begin{theorem}
    \href{https://github.com/Thmoas-Guan/mathlib4_fork/blob/ABS-Criterion-Project-new/Mathlib/RingTheory/RegularLocalRing/RegularRing/GlobalDimension.lean#L25-L64}{\faExternalLink}
    For a regular ring $R$, $\gldim(R) = \dim(R)$.
\end{theorem}
\begin{proof}
This is a direct consequence of Lemma ~\ref{lemma:gldimeqsup} and Theorem ~\ref{thm:reggldim}.
\end{proof}

Since localization of a regular local ring at a prime is regular, regularity of a ring can be verified on maximal ideals only.
\begin{theorem}
    \href{https://github.com/Thmoas-Guan/mathlib4_fork/blob/ABS-Criterion-Project-new/Mathlib/RingTheory/RegularLocalRing/RegularRing/Localization.lean#L41-L67}{\faExternalLink}
    A ring $R$ is regular if and only if $R_\mathfrak{m}$ is a regular local ring for every maximal ideal $\mathfrak{m}$.
\end{theorem}

We are finally able to prove that a regular local ring is a regular ring.\href{https://github.com/Thmoas-Guan/mathlib4_fork/blob/ABS-Criterion-Project-new/Mathlib/RingTheory/RegularLocalRing/RegularRing/Localization.lean#L36-L39}{\faExternalLink}

\subsubsection{\textbf{Hilbert's Syzygy Theorem}}\label{subsec:Syzygy}

With all the previous results in this section we are able to prove a weak version of Hilbert's Syzygy Theorem.
\href{https://github.com/Thmoas-Guan/mathlib4_fork/blob/ABS-Criterion-Project-new/Mathlib/RingTheory/RegularLocalRing/RegularRing/Syzygy.lean#L76-L81}{\faExternalLink}
\begin{theorem}[Hilbert's Syzygy,{\cite[Corollary 18.1]{mats_commalg_2ed}}]
    For a field $k$, $\gldim(k[x_1, \cdots, x_n]) = n$.
\end{theorem}
\begin{lstlisting}
theorem Hilberts_Syzygy (k : Type u) [Field k] [Small.{v} k] (n : ℕ) :
    globalDimension.{v} (MvPolynomial (Fin n) k) = n
\end{lstlisting}
However, due to the current setup, we still have some distance from the original theorem first stated in \cite{Hilbert1890}.
To prove the weak version above, it only remains to show that a polynomial ring over a regular ring is regular.
\begin{theorem}
    \href{https://github.com/Thmoas-Guan/mathlib4_fork/blob/ABS-Criterion-Project-new/Mathlib/RingTheory/RegularLocalRing/RegularRing/Polynomial.lean#L90-L163}{\faExternalLink}
    For a regular ring $R$, $R[X]$ is regular.
\end{theorem}
Because both regularity and the Cohen--Macaulay property are local conditions, the proof reduces similarly to the case over a regular local ring $(R, \mathfrak{m}, k)$ and localizing at a prime $\mathfrak{p}$ of $R[X]$ containing $\mathfrak{m}$. We managed to find some common lemmas that can be reused through the process.\href{https://leanprover-community.github.io/mathlib4_docs/Mathlib/RingTheory/Ideal/MonicSpan.html#Polynomial.exists_monic_span_sup_map_eq}{\faExternalLink}

\section{Discussion}\label{sec:discuss}

\subsection{Remarks About Universe Levels}\label{subsec:design}

In the formalization of category theory in \Leanf{}, one often needs to carefully consider universe levels.
However, readers interested only in the mathematical content may assume throughout that the ring and the modules are all in the same universe $u$, which is the most important case from the mathematical point of view. 
In this situation, the $\Ext$ groups also live in universe $u$, and no comparison between different universe levels is needed.

But there might be situations that do not arise in this way. For example, this issue arises when proving that for an arbitrary field $K$ of positive characteristic $p$, there exists a Cohen ring whose residue field is isomorphic to $K$. For a local ring $(R, \mathfrak{m}, k)$ and $L$ a field extension of $k$, we have a lemma for obtaining a flat local $R$-algebra with residue field $L$. It is proven by category theory, hence we restrict $R$ and $L$ to the same universe in the first place. We want to apply it to $\mathbb{Z}_p$ and $K$, but $\mathbb{Z}_p$ lives in \lean{Type} while we wish $K$ to be in an arbitrary universe level. In this case, we have to forcefully generalize the universe level in the lemma to make it applicable in this case.

Due to the existence of the situations above, in this paper, we work in a slightly more general setting to make our results directly applicable. We allow \lean{R : Type u} and work with modules in \lean{ModuleCat.\{v\} R}, under the assumption \lean{[Small.\{v\} R]}. This assumption ensures that \lean{ModuleCat.\{v\} R} has enough projectives, allowing us to consider the $\Ext$ groups in universe $v$. 
When $R$ or $R/I$ is regarded as an object of \lean{ModuleCat.\{v\} R}, we use \lean{Shrink} to replace it by an isomorphic representative whose carrier lies in universe $v$. 
Thus, the occurrences of \lean{Small} and \lean{Shrink} are devices for universe bookkeeping and are not part of the mathematical content of the results.

However, as we try not to restrict everything--especially every ring to the same universe level, it does require some effort to prove some homological invariants are stable across universe levels, for example $\gldim$ in Subsection \ref{subsec:gldimpre}. This is why we usually prefer to make definitions in non-categorical language. However, in the case of depth, as we see in Theorem \ref{thm:Ischbeck}, we are indeed benefiting from the general definition of depth, we still choose to define it using $\Ext$. But as we usually work over a Noetherian ring and finitely generated modules, we have the regular sequence characterization, in this case it just wouldn't cause much inconvenience.

\subsection{Related Work}\label{subsec:relwork}

As far as we are aware, in other theorem provers the concepts of derived categories and derived functors are often not formalized, let alone theories built on them such as $\depth$ and homological dimension. In \Leanf{}, this work also represents a first step toward applications of homological methods in commutative algebra.

Before this work, formalizations of commutative algebra in \Leanf{} were primarily motivated by applications in algebraic geometry and algebraic number theory. For instance, the formalization of homogeneous ideals originated from the formalization of projective space \cite{proj2023}; developments on Dedekind domains and class groups were motivated by their roles in algebraic number theory \cite{classgroup}; and it is the formalization of local fields that furthered the formalization of valuation rings \cite{localfield}. In contrast, the present work undertakes a systematic formalization of commutative algebra itself in combination with homological algebra, positioning this development as a foundational body of theory on which algebraic geometry and algebraic number theory can subsequently rely.

Among formalizations of commutative and homological algebra in \Leanf{} directly related to this project, there is a project led by Christian Merten that formalizes the fact that the Krull dimension of a polynomial ring over a Noetherian ring increases by exactly one \href{https://github.com/leanprover-community/mathlib4/pull/27542}{\faExternalLink}. That proof proceeds by showing (a version of) \cite[Theorem 19]{mats_commalg_2ed}\href{https://github.com/leanprover-community/mathlib4/pull/27510}{\faExternalLink}. We would like to thank Andrew Yang here for his efforts in commutative algebra, scattered throughout \mathlib{}, which greatly benefited this work.

Regarding homological algebra, the Liquid Tensor Experiment \cite{lean2022_LTE} motivated and relied on early theories of derived functors defined via resolutions. The $\Ext$ that we mainly use in this article was developed later, along with theory of derived category \cite{derived_cat}.

Even though we used almost only the definition of $\Ext$ and its long exact sequence, the work done at the level of derived categories has already laid a solid foundation for the subsequent formalization.

\subsection{Future Work}

Building on the development of Cohen--Macaulay rings and regular local rings, there are several natural directions for future work.

\paragraph{\textbf{Cohen--Macaulay Modules over Arbitrary Rings}}
For a non-local ring $R$, one can define Cohen--Macaulay modules as follows: an $R$-module $M$ is Cohen--Macaulay if for every prime ideal $\mathfrak{p}$ of $R$, the localized module $M_{\mathfrak{p}}$ is Cohen--Macaulay over $R_{\mathfrak{p}}$.

\paragraph{\textbf{Miracle Flatness}} The following is known as miracle flatness.
\begin{theorem}[{\cite[00R4]{stacks-project}}]
For a local homomorphism $R \to S$ with $R$ a regular local ring and $S$ a Cohen--Macaulay local ring, if $\dim S = \dim R + \dim(S/\mathfrak{m}_R S)$, then $R \to S$ is flat.
\end{theorem}
After the developments in this project, most of the required techniques are already formalized or are only a small step away. The remaining obstacle is formalizing a variant of the local criterion for flatness:
\begin{lemma}[{\cite[00ML]{stacks-project}}]
Let $R \rightarrow S$ be a local homomorphism of Noetherian local rings. Let $I \neq R$ be an ideal in $R$. Let M be a finitely generated $S$-module. If $\Tor_R^1(M,R/I)=0$ and $M/IM$ is flat over $R/I$, then $M$ is flat over $R$.
\end{lemma}

\paragraph{\textbf{Serre's Criterion of Normality}} We use the following standard notation for a ring $R$:
\begin{itemize}
    \item ($S_k$) : $\depth(R_\mathfrak{p}) \geq min\{k, \height(\mathfrak{p})\}$ for all $\mathfrak{p} \in \Spec(R)$,
    \item ($R_k$) : if $\mathfrak{p} \in \Spec(R)$ and $\height(\mathfrak{p}) \leq k$, $R_\mathfrak{p}$ is regular.
\end{itemize}
Given our formalization of $\depth$, we can state the criterion of normality as follows: for a Noetherian ring $R$, $R$ is normal if and only if it satisfies ($S_2$) and ($R_1$) (see \cite[Theorem 39]{mats_commalg_2ed}). Formalizing this criterion would be an interesting project that will develop many results about normal domains along the way.

\paragraph{\textbf{Other Versions of Hilbert's Syzygy Theorem}}

In Subsection ~\ref{subsec:Syzygy}, we proved that the global dimension of $k[x_1,\dots,x_n]$ is $n$, which is only slightly stronger than the statement that every finitely generated module over $k[x_1,\dots,x_n]$ has a projective resolution of length $n$—the weakest version of Hilbert's Syzygy Theorem. These projective resolutions can indeed be upgraded to free resolutions using the Quillen--Suslin theorem (see \cite{Quillen1976}\cite{Suslin1976}) or other constructions specific to polynomial rings over fields (see \cite[Corollary 19.8]{Eisenbud2022}). Additionally, there is a graded version of this theorem (see \cite[Corollary 19.7]{Eisenbud2022}); however, with the current setup, we still have some distance to the formalization of the graded version.

\bibliographystyle{amsalpha}
\bibliography{main}

@misc{stacks-project,
  author       = {The {Stacks project authors}},
  title        = {The Stacks project},
  howpublished = {\url{https://stacks.math.columbia.edu}},
  year         = {2025},
}

@book{mats_commalg_2ed,
  author    = {Hideyuki Matsumura},
  title     = {Commutative Algebra},
  edition   = {Second},
  note      = {Revised and modernized edition by {\TeX}romancers; original edition 1980; this edition 2022},
  year      = {2022},
}

@book{mats_commring,
author = {Hideyuki Matsumura},
title = {Commutative Ring Theory},
translator = {Miles Reid},
series = {Cambridge Studies in Advanced Mathematics},
number = {8},
publisher = {Cambridge University Press},
year = {1987},
note = {Print publication year 1987; online publication June 2012},
isbn = {9781139171762},
doi = {10.1017/CBO9781139171762},
url = {https://doi.org/10.1017/CBO9781139171762}
}

@book{CM_ring,
  author    = {Winfried Bruns and J\"urgen Herzog},
  title     = {Cohen--Macaulay Rings},
  series    = {Cambridge Studies in Advanced Mathematics},
  number    = {39},
  edition   = {Revised},
  publisher = {Cambridge University Press},
  address   = {Cambridge, UK},
  year      = {1998},
  isbn      = {0-521-56674-6},
  note      = {First paperback edition with revisions; original edition 1993},
}

@book{Eisenbud2022,
  author    = {Eisenbud, David},
  title     = {Commutative Algebra: With a View Toward Algebraic Geometry, Second Edition},
  publisher = {Springer},
  year      = {2022}
}

@misc{haution_homalg,
  author       = {Olivier Haution},
  title        = {Homological Methods in Commutative Algebra},
  howpublished = {Lecture notes, Ludwig-Maximilians-Universit{\"a}t M{\"u}nchen},
  month        = {Summersemester},
  year         = {2017},
  note         = {Available at: \url{https://haution.gitlab.io/pdf/Local_Algebra.pdf}},
  url          = {https://haution.gitlab.io/pdf/Local_Algebra.pdf},
}

@book{Macaulay1916,
  author    = {Macaulay, F. S.},
  title     = {The Algebraic Theory of Modular Systems},
  publisher = {Cambridge University Press},
  year      = {1916},
  address   = {Cambridge}
}

@article{NorthcottRees1954,
  author    = {Northcott, D. G. and Rees, D.},
  title     = {Reductions of ideals in local rings},
  journal   = {Mathematical Proceedings of the Cambridge Philosophical Society},
  year      = {1954},
  volume    = {50},
  pages     = {145--158},
  doi       = {10.1017/S0305004100020952}
}

@article{Hilbert1890,
  author    = {Hilbert, David},
  title     = {Ueber die Theorie der algebraischen Formen},
  journal   = {Mathematische Annalen},
  year      = {1890},
  volume    = {36},
  pages     = {473--534},
  doi       = {10.1007/BF01208503}
}

@inproceedings{Serre1956,
  author    = {Serre, Jean-Pierre},
  title     = {Sur la dimension homologique des anneaux et des modules noeth\'eriens},
  booktitle = {Proceedings of the International Symposium on Algebraic Number Theory, Tokyo {\&} Nikko, 1955},
  pages     = {175--189},
  publisher = {Science Council of Japan},
  address   = {Tokyo},
  year      = {1956},
  language  = {French}
}

@article{AuslanderBuchsbaum1956,
  author    = {Auslander, Maurice and Buchsbaum, David A.},
  title     = {Homological dimension in Noetherian rings},
  journal   = {Proceedings of the National Academy of Sciences of the United States of America},
  year      = {1956},
  volume    = {42},
  pages     = {36--38},
  doi       = {10.1073/pnas.42.1.36}
}

@article{AuslanderBuchsbaum1957,
  author    = {Auslander, Maurice and Buchsbaum, David A.},
  title     = {Homological Dimension in Local Rings},
  journal   = {Transactions of the American Mathematical Society},
  year      = {1957},
  volume    = {85},
  pages     = {390--405},
  doi       = {10.2307/1992962}
}

@article{Ischebeck1969,
  author  = {Friedrich Ischebeck},
  title   = {Eine Dualit{\"a}t zwischen den Funktoren {Ext} und {Tor}},
  journal = {Journal of Algebra},
  year    = {1969},
  volume  = {12},
  pages   = {510--531},
  doi     = {10.1016/0021-8693(69)90090-8}
}

@article{Quillen1976,
  author    = {Quillen, Daniel},
  title     = {Projective Modules over Polynomial Rings},
  journal   = {Inventiones Mathematicae},
  year      = {1976},
  volume    = {36},
  number    = {1},
  pages     = {167--171},
  doi       = {10.1007/BF01403183}
}

@article{Suslin1976,
  author    = {Suslin, Andrei A.},
  title     = {Projective Modules over Polynomial Rings are Free},
  journal   = {Doklady Akademii Nauk SSSR},
  year      = {1976},
  volume    = {229},
  number    = {5},
  pages     = {1063--1066},
  note      = {English translation: \emph{Soviet Math. Dokl.} \textbf{17} (1976), 1160--1164}
}

@article{intersection65,
author = {Jean-Pierre Serre},
title = {Alg`ebre locale. Multiplicit'es},
journal = {Lecture Notes in Mathematics},
volume = {11},
publisher = {Springer-Verlag},
year = {1965},
}

@article{deformation87,
author = {Murray Gerstenhaber and Samuel D. Schack},
title = {Algebras, bialgebras, quantum groups, and algebraic deformations},
journal = {Mathematics and Its Applications},
volume = {10},
pages = {9--35},
year = {1987},
}

@article{Auslander1953,
  author    = {Maurice Auslander},
  title     = {On the dimension of modules and algebras. I},
  journal   = {Nagoya Mathematical Journal},
  volume    = {2},
  year      = {1953},
  pages     = {64--76},
  issn      = {0027-7630}
}

@inproceedings{mathlib,
author = {The {Mathlib Community}},
editor={Blanchette, Jasmin C. and Hrițcu, Catalin},
title = {The Lean Mathematical Library},
year = {2020},
isbn = {9781450370974},
publisher = {ACM},
doi = {10.1145/3372885.3373824},
booktitle = {CPP 2020},
pages = {367–381},
numpages = {15},
}

@misc{lean2022_LTE,
  author       = {Johan Commelin and Adam Topaz and others},
  title        = {Completion of the Liquid Tensor Experiment},
  howpublished = {Blog post},
  year         = {2022},
  url          = {https://leanprover-community.github.io/blog/posts/lte-final/},
  note         = {Lean community blog.},
}

@article{derived_cat,
  author  = {Jo\"el Riou},
  title   = {Formalization of derived categories in Lean/mathlib},
  journal = {Annals of Formalized Mathematics},
  volume  = {1},
  pages   = {1-42},
  year    = {2025},
  month   = {jul},
  day     = {15},
  doi     = {10.46298/afm.13609},
  url     = {https://doi.org/10.46298/afm.13609},
  note    = {Source: HAL:hal-04546712v5. Submitted: 2024-05-17; Accepted: 2025-02-11. License: CC BY 4.0},
}

@inproceedings{proj2023,
author = {Jujian Zhang},
title = {Formalising the Proj Construction in Lean},
booktitle = {14th International Conference on Interactive Theorem Proving (ITP 2023)},
year = {2023},
publisher = {Leibniz International Proceedings in Informatics (LIPIcs)},
address = {Dagstuhl, Germany},
doi = {10.4230/LIPIcs.ITP.2023.35},
note = {Department of Mathematics, Imperial College London, UK; Supplementary software: https://github.com/leanprover-community/mathlib/pull/18138/commits/00c4b0918a2c7a8b62291581b0e1eddf2357b5be},
}

@article{classgroup,
author = {Anne Baanen and Sander R. Dahmen and Ashvni Narayanan and Filippo A. E. Nuccio Mortarino Majno di Capriglio},
title = {A Formalization of Dedekind Domains and Class Groups of Global Fields},
journal = {Journal of Automated Reasoning},
year = {2022},
volume = {66},
pages = {611--637},
doi = {10.1007/s10817-022-09644-0},
url = {https://doi.org/10.1007/s10817-022-09644-0},
}

@inproceedings{localfield,
author = {María In{'e}s de Frutos-Fern{'a}ndez and Filippo Alberto Edoardo Nuccio Mortarino Majno di Capriglio},
title = {A Formalization of Complete Discrete Valuation Rings and Local Fields},
booktitle = {Proceedings of the 13th ACM SIGPLAN International Conference on Certified Programs and Proofs (CPP'24)},
year = {2024},
month = {January},
address = {London, UK},
publisher = {ACM},
doi = {10.1145/3636501.3636942},
url = {https://doi.org/10.1145/3636501.3636942},
eprint = {2310.01998},
eprinttype = {arXiv},
eprintclass = {cs.LO},
note = {Also available as arXiv:2310.01998, v2 (08 Dec 2023)},
}

\end{document}